\documentclass{amsart}

\usepackage{amsmath, amsthm}
\usepackage{amssymb, latexsym}
\usepackage{amsfonts}

\newtheorem{theorem}[equation]{Theorem}

\newtheorem{proposition}[equation]{Proposition}
\newtheorem{lemma}[equation]{Lemma}

\newtheorem{corollary}[equation]{Corollary}

\newtheorem{Remark}[equation]{Remark}
\newtheorem{remark}[equation]{Remark}

\newtheorem{convention}[equation]{Convention}

\newcounter{com}

\newtheorem{comments}{}

\newcommand{\beql}[1]{\begin{equation}\label{#1}}
\newcommand{\eeq} {\end{equation}}

\font\Aaa=msam10
\def\nor{\hbox{ {\Aaa  C} }}

\def\qed{\hbox{~~\Aaa\char'003}}

\font\Bbb=msbm10
\def\semi{\hbox{\Bbb o}}

\def\Z{\hbox{\Bbb Z}}

\def\F{\hbox{\Bbb F}}

\numberwithin{equation}{section}

\let\define=\def
\let\redefine=\def

\def\Tr{{\rm Tr}}
\def\Aut{{\rm Aut}}

\define\a{{ \alpha }}
\redefine\b{{ \beta }}

\redefine\d{ \delta }
\define\e{{ \epsilon }}
\redefine\g{\gamma}
\def\G{\Gamma}

\redefine\D{{ \Delta }}
\redefine\l{{ \lambda }}

\define\({ \left( }
\define\){ \right) }
\define\[{ \left[ }
\define\]{ \right] }
\define\<{ \langle }
\define\>{ \rangle }
\let\ljunk=\{
\let\rjunk=\}
\redefine\{{\left\ljunk}
\redefine\}{\right\rjunk}

\redefine\.{ \diamomd }
\redefine\+{\oplus}
\redefine\.{ \cdot}
\redefine\#{\sharp}

\redefine\.{ \cdot }
\redefine\l{ \lambda }

        \def\PSL{{\rm PSL}}
        \def\SL{{\rm SL}}

        \def\SU{{\rm SU}}

        \def\Sp{{\rm Sp}}
        
        \def\GL{{\rm GL}}
      \def\Spin{{\rm Spin}}

     \def\diag{{\rm diag}}

        \def\Sz{{\rm Sz}}

        \font\Aaa=msam10
        \def\nor{\hbox{\Aaa C}}

\font\Aaa=msam10
\def\nor{\hbox{ {\Aaa  C} }}

\def\qed{\hbox{~~\Aaa\char'003}}

\font\Bbb=msbm10
\def\semi{\hbox{\Bbb o}}

\def\Z{\hbox{\Bbb Z}}

\def\F{\hbox{\Bbb F}}

\def\a{\alpha}
\def\b{\beta}

\def\div{ \kern-.5pt\hbox{\big |} }
\def\ndiv{ {\not\kern-.5pt\hbox{\big |}\,} }
\def\ndivv{ {\not\kern+1.5pt\hbox{$\mid$}\,} }

\def\lr#1{{\langle #1 \rangle}}

\def\Tr{{\rm T}}

\def\Aut{{\rm Aut}}

\def\GL{{\rm GL}}
\def\norm{{\rm N}}
\def\cent{{\rm C}}
\def\Oh{{\rm O}}
\def\zent{{\rm Z}}
\def\diag{{\rm diag}}

\def\ll{{\varpi}}

\def\B{^2\kern-.8pt B}
\def\G{^2\kern-.8pt G}
\def\EH{^2\kern-.8pt\hat  E}
\def\E{^2\kern-.8pt E}
\def\D{^3\kern-1pt D}
\def\FF{^2\kern-.8pt F}

\def\rel{\vskip.01in\medskip \noindent{\bf Reliability:}~}

\newdimen\refcodesize
\newbox\seriesbox
\def\para#1{\medskip\noindent{\bf #1}~}

\def\timing{\medskip\noindent{\bf Time:}~}

\def\proof{\noindent {\bf Proof.~}}

\def\PSL  {{\rm PSL }}

\def\Po{{\rm P\Omega}}

\def\Sp{{\rm Sp}}

\def\GL{{\rm GL}}
\def\SU{{\rm SU}}

\def\SL{{\rm SL }}
\def\GL{{\rm GL}}
\def\Spin{{\rm Spin}}

\def\SSS{{\mathcal S}}

\def\ad  {{\rm ad }}

\def\ppd{{\rm ppd}}

\begin{document}

\title[Black box exceptional  groups of Lie type]
{Black box exceptional  groups of Lie type}

\thanks{This research was supported in part by
 NSF grants  DMS~9731421, DMS~0242983 and DMS~0753640,
 and 
NSA grant MDA-9049810020.}

       \author{W. M. Kantor}
       \address{University of Oregon,
       Eugene, OR 97403}
       \email{kantor@uoregon.edu}

    \author{K. Magaard}
       \address{University of Birmingham,  Edgbaston, Birmingham B15 2TT}
       \email{k.magaard@bham.ac.uk}

\subjclass[2000]{Primary: 20D06, 20G40 ~
Secondary: 20B40, 20G41, 20P05, 68Q25}

{\abstract 
11111
}

\begin{abstract}
If a black box  group   is known
to be isomorphic to an exceptional simple group of Lie type 
of (twisted) rank $>1$, other than any  $\FF_4(q)$,
over a field
of known size, a Las Vegas algorithm is given to produce a constructive
isomorphism. In view  of its timing, this algorithm yields an upgrade of all known nearly linear time Monte
Carlo permutation group algorithms to Las Vegas algorithms when the input
group has no composition factor isomorphic to  
 any group $\FF_4(q)$ or $^2G_2(q)$.

\vspace{-22pt}

 \end{abstract}
  
\maketitle

\section {Introduction}
\label {Introduction}  

In a number of algorithmic settings it is essential to
take a permutation group or matrix group that is known (probably) to be simple and
produce an explicit isomorphism with an explicitly defined simple
group, such as a group of matrices 
  (\cite{LG,KS1,Ka3} contain background on this and many related questions).
This has been accomplished for the much more general setting of black box
classical groups in
\cite{KS1,Br1,Br2,Br3,BrK1,BrK2,LMO} (starting with the  groups
$\PSL(d,2)$  in
\cite{CFL}). Black box alternating groups are dealt with in
\cite{BLNPS}.  In this paper we consider this  identification question
for black box exceptional groups of Lie type.
Note that the {\em name} of the group can be found quickly using Monte Carlo algorithms in suitable settings
\cite{BKPS,KS4,KS5,LO}.

The elements of a black box group $G$ are assumed to be encoded
by 0-1 strings of uniform length, and $G$ is specified as
$G=\lr\SSS$ for some set $\SSS$ of elements of $G$. 
%\newpage
Our main  result is as follows
(where $\epsilon $  is 1 in general%,  
~and 3 for 
$\D_4(q)$):

\begin{theorem}
\label{Main Theorem} \label{main theorem}
There is a Las Vegas algorithm which$,$ when given 
a  black box group $G=\lr\SSS$ 
isomorphic to a
perfect central extension of a
 simple exceptional group of Lie type
of $($twisted$)$  rank $>1$ and given field size $q,$ other than 
any 
$\FF_4(q),$
 finds the following$:$ 

\begin{itemize}

\item[\rm(i)] The name of the simple group  of Lie type to which
$G/\zent(G)$ is isomorphic$;$  and  	

\item[\rm(ii)] A new set $\SSS ^*$ generating $G,$
 a generating set $\hat \SSS$ of the universal cover $\hat G$ of the
simple group in {\rm (i)}
and an epimorphism
$\Psi\colon\! {\hat G}\to G,$  specified by the requirement that  $\hat
\SSS\Psi =\SSS ^*$.
\end{itemize} 
Moreover$,$  the data structures underlying {\rm(ii)} yield 
algorithms for each of the following$:$ 
\begin{itemize} 

\item[\rm(iii)] Given $g\in G,$  find
$\hat g\in\hat   G $ such that 
$g= \hat g\Psi ,$ and a straight-line 
program  of length%
\footnote{All logarithms are to the base $2$.}
 $O(\log q)$ from
$  \SSS ^*$ to  $g;$  and    

\item[\rm(iv)]   Given $\hat g\in {\hat G},$ find $\hat g\Psi $ and a
straight-line program of length $O(\log q)$ 
from $ \hat \SSS  $ to $\hat g.$%

\end{itemize}   
In addition$,$  the following all hold.
\begin{itemize}

\item[\rm(v)] $\SSS^*$ has size   $O(\log q)$
 and contains a generating
set for $G$ consisting of  root    elements.  

\item[\rm(vi)] The  algorithm for {\rm(ii)} is an
$O(\xi q ^\epsilon \log q+\mu q ^\epsilon  \log ^ 2\hspace{-1pt}q   ) $-time
Las Vegas algorithm   succeeding with probability  $>1/2,$  where
$\mu $   is an upper bound on the time required for
each group operation in $G,$ and
$\xi\ge\mu $ is an upper bound on the time requirement per element for the
construction of independent$,$ $($nearly$)$ uniformly distributed  
random elements of
$G.$ 

In additional $O \big(|\SSS|\log |\SSS| (\xi q^\epsilon \log q+\mu q^\epsilon \log ^ 2\hspace{-1pt}q ) \big)$ time 
it can be verified that $G$ is indeed
isomorphic to a perfect  central extension of the  
exceptional group in~{\rm (i)}.

\item[\rm(vii)] 
 The algorithm  
for  {\rm(iii)} is   Las Vegas$,$  running in 
$O(\xi q^\epsilon \log q+\mu q^\epsilon \log ^ 2\hspace{-1pt}q   )$ time and 
 succeeding with probability  $>1/2;$  while the  algorithm for 
 {\rm(iv)} is   deterministic and runs in  
$O( \mu\log q)$  time.

\item[\rm(viii)] The center of $G$ can be found in $O( \mu  \log q )$
time.
\end{itemize} 
\end{theorem} 
% \noindent
Parts (ii-iv) are the requirements  for a 
\emph{constructive  epimorphism}  $\Psi\colon\! {\hat G}\to G$. 
The verification at the end of (vi)  is omitted in
some references, since
$G$ is {\em assumed} to be an epimorphic image of a specific group $\hat
G$ which, in turn, is isomorphic to (a central extension of) an
explicitly constructed  subgroup $G_0$ of $G$ (cf.
Proposition~\ref{constructing lambda}).   In practice, it is hard to
imagine that this test would be omitted since it appears to be the only
way to guarantee that the group $G$ behaves  as hypothesized. 
We note that, in (iv),  $\hat g\in \hat G$ might be given  in standard Bruhat normal form but alternatively  might merely be given as an automorphism of the associated Lie algebra (cf. Remark~\ref{elements of hat G}).
 It is also worth remarking that we use $\e=1$ for groups of type $\E_6(q)$,
since that is the case for its Levi factor $\SU(6,q)$ 
(Theorem~\ref{previous algorithms}).

The above algorithms do not run in polynomial time:   the timing in (vi) and (vii) have factors $q$.  At present there 
are no polynomial-time algorithms for the type of problem considered here, neither in the black box setting nor even 
in the matrix group one. 
This was already evident for classical groups in \cite{KS1} and,  even earlier, in \cite{CLG2}. A standard way around this obstacle involves a lovely idea in \cite{CoLG} (used in 
\cite{BrK1,BrK2,Br2,Br3,LMO}):  \emph{use suitable oracles}.  The preceding references {\em assume}  the availability of  
 an oracle  that constructively 
recognizes subgroups $\SL(2,q)$.  
This was motivated by \cite{CoLG}, which  deals with matrix groups 
 and  assumes  the availability of a  {\em
Discrete Log  oracle} for $\F_q^*$.  In this matrix group setting,   
\cite{CoLG,LGO}  provide  a constructive Las Vegas
algorithm   for a group isomorphic to  a nontrivial homomorphic image of $\SL(2,q)$  in any irreducible
representation in characteristic dividing $q$, running in time
that is   polynomial in the input length,
\emph{assuming} the availability of a Discrete Log oracle. 
In effect,  this idea  replaces annoying factors  $q$ by an  oracle.
This is discussed further in Section~\ref{Concluding remarks}, 
Remarks~\ref {The factor $q$.}--\ref{Involution centralizers.}, making it 
 clear that this will not be the last paper on this type of problem!

 A rough outline of the proof of the theorem  is given in Section~\ref{Outline}. 
The first part   resembles 
\cite{KS1}: we find a long root element, then build  a subgroup
$\SL(3,q)$, and also a subgroup 
  $\Spin_8^{-}(q)$ when the Lie  rank is more than 2.  We then use
pieces of  these groups to obtain   the centralizer of a subgroup
$\SL(2,q)$ generated by long root groups. However, there is no 
module to aim for that is as nice as in the classical case.  Hence,
instead we proceed directly to obtain all of the root groups
corresponding to a root system, and then verify the standard commutator
relations that  define these groups: the corresonding  presentation guarantees the
Las Vegas nature of our algorithm.   

Our proofs are divided into two parts, with  rank 
$>2$ and rank  2     in Sections~\ref{Rank
$>2$} and 
\ref{Rank $2$ groups}, respectively.
Section~\ref{Concluding remarks}   contains
remarks concerning improvements or variations on the theorem and
the algorithms.  

In view of  \cite{KS2} (and \cite{BrB}), we obtain the following immediate but
significant consequence of the above theorem:

\begin{corollary}
\label {Upgrade corollary}
Given a permutation group $G \le S_n$ with no composition factor
isomorphic to    
 any group $\FF_4(q)$ or $^2G_2(q),$ all known nearly 
linear time Monte Carlo algorithms dealing with $G$ can be upgraded to
Las Vegas algorithms.   
\end{corollary}
 
The stated algorithms find $|G|$ and a composition series of $G$, 
among many other things (cf. \cite{Ser}).
In fact, it can be shown that the groups $^2G_2(q)$ do not need to be excluded here;
see  Section~\ref{Concluding remarks}, 
Remark~\ref{Rank 1 groups.}.

\subsection{Background}
\label{Background}  
For background on groups of Lie type we refer 
to   \cite{Ca1,GLS}.
For background on  required aspects of black box groups, 
in particular  for discussions of the parameters $\xi$ and $\mu $ in the
theorem, see
\cite[{\bf2.2.2}]{KS1}.  Thus,     
{\em we assume that} $\xi \ge \mu |\SSS|$ and $\mu\ge {\tt N} $
if ${\tt N}$ is the string length of the elements of our black box group $G$.
Moreover, 
${\tt N} \ge \log |G| >C \log q$ for some constant $C$, 
since we are dealing with  exceptional groups of Lie type over 
$\F_q.$

 We note that, as in
\cite[{\bf2.2.4}]{KS1}, we presuppose the availability of  independent  $($nearly$)$
uniformly distributed random elements of $G$, a major result in  \cite{Bab}
(compare \cite{CLMNO,Dix}). 

 {\em Straight-line programs} from $\SSS$ to  elements of $ G = \langle \SSS \rangle$ are also defined and discussed in
\cite[{\bf2.2.5}]{KS1}.   For use in \cite{KS2}
(or   in Corollary~\ref{Upgrade corollary}), part (iii) of
the theorem    needs the stated straight-line program, not just the preimage
$\hat g $.

In general the symbol $\ppd^\#(p;n)$ stands for some integer divisible
by a prime~$r$ (a \emph{primitive prime divisor}) such that $r\div p^{n}-1$ but  $r\ndiv p^{i}-1$ for $1\le i<n$ (cf.
\cite{Zs}).
The exceptions to~this definition are:~$\ppd^\#(p;1)$ with $p>5$ a Fermat prime, where we require  
divisibility by 4; 
$\,\ppd^\#(p;2)$ with $p$ a Mersenne prime, where we require  
divisibility by 4;   and 
$\ppd^\#(2;6)$,  where we require
divisibility by 21.
It is easy to test this property of an integer for a single $\ppd^\#(p;n)$ requirement
\vspace{1pt}
(\cite[p.~578]{NP},
\cite[Lemma~2.7]{KS1}), and hence also for a product
$\ppd^\#(p;n_1)\cdots\ppd^\#(p;n_k)$  of a bounded number  of them
(where $n_1<\cdots<n_k$). 
In those references, the time requirement for such tests is far smaller than other aspects of
our algorithms, and hence will be ignored.

\para{
Notation:}  We always write $q=p^e$, where $p$  is  the characteristic of $G$.

We will usually have available a field $\F= \F _q$ obtained from subgroups
of $G$; and also an extension field $\F'$ of $\F $  of degree 1, 2 or 3.
We choose an $\F_p$--basis $\{f_1,\ldots\}$ of $\F'$ such that   $f_1=1$ and 
 $\{f_1,\ldots,f_e\}$ is a basis of $\F$.
 
 In view of the discussion in \cite[Sec.~2.3]{KS1}, we will always assume that {\em field operations  can be carried out in constant
time.}

\subsection{Outline}
\label{Outline of algorithm}
\label{Outline}

A very rough summary of our approach to 
Theorem~\ref{main theorem}(ii,vi)  is as follows (with many details suppressed
or ignored).

%\vspace{-6pt}
\begin {itemize}
\item
Use random group elements and primitive prime divisors to find
$\tau\in G$ of special order, in particular such that some power
$z=\tau^l$ is a long root element
(cf. Sections~\ref {primitive prime divisors}, \ref{Finding a root element},
\ref{Root element and S}).  (In types $E_7$ and $E_8$ we need two such elements
$\tau$ of different specific orders.)
%\vspace{-6pt}

\item
Find three conjugates of $z$  that generate a subgroup
$S=\SL(3,q)$ (cf. Sections~\ref{Probability},
\ref{Finding a root group}, \ref{Root element and S}),
together with a    subgroup $ R\cong \SL(2,q)$ of $S$   also
generated by conjugates of $z$.
Much of the algorithm depends heavily on $\SL(2,q)$ and $\SL(3,q)$
subgroups.
%\vspace{-6pt}

\item 
For rank  $>2$  use $S$ and a conjugate of $z$ to construct a 
$G$-conjugate of $R$ lying in
$L=\cent_G(R)$ (cf. Section~\ref{The long subgroups  $J$ and    $R_1$.} );  this $\SL(2,q)$ and variants of  the element(s) 
$\tau$ are used to generate
$L$ (cf. Section~\ref{Finding  $L$}).  If the rank is  $2$  then  $\cent_S(R)$   and $\tau$ generate $L$
(which this time is an  $ \SL(2,q^\e)$; cf. Section~\ref{The subgroups $L$ and $Q$}).

This heavily depends on  the uniqueness of the triple $(R,S,L  )$ up to conjugacy in $G$.

%\vspace{-6pt}

\item
Find a (maximally) split torus $T$ normalizing $L$ and $S$.
 Use it to 
construct root systems of
$L$ and $S$ with respect to the   tori $T\cap L$ and $T\cap S$.
Use  commutators of  root groups of $S$ and $L$ to 
find    root groups and a root system $\Phi_G$ for $G$
(cf. Sections~\ref {Matching},
\ref{The subgroups $L$ and $Q$}).

%\vspace{-6pt}
\item
The new generating set $\SSS^*$ for $G$ contains the union of sets of
generators of these root groups $X_\alpha$,  $\alpha \in \Phi_G$.
Verify a version of the Steinberg presentation \cite{St} for the subgroup $G_0$
generated by these subgroups $X_\alpha$
(cf. Sections~\ref{Finding and verifying the presentation}, 
\ref{Root groups.}).  
%\vspace{-6pt}
\item
Show that each of the given generators for $G$ is in
$G_0$, so that $G_0=G$ (cf. Sections~\ref {Straight-line programs},
\ref{Proof of Theorem}).

\end {itemize}
\vspace{-3pt}
%\newpage

\subsection{Recognition algorithms used}
\label{Recognition algorithms used}
We will use existing algorithms for constructive recognition of various black box groups.  Since their timing is crucial for us, we state the instances and timings in the next result, which refers to the 
counterparts in our Theorem~\ref{main theorem}: 

\begin{theorem}
\label{previous algorithms}
Let $G=\langle \SSS \rangle$ be a black box group that 
is isomorphic to a nontrivial homomorphic image of 
$\SL(2,q),$
$\SL(3,q),$ $ \Sp(6,q),$ $ \SU(6,q),$
$\Spin^-_8(q)$ or   $\Spin^+_{12}(q)$.  Then there are  
algorithms for  the natural analogues of {\:\rm Theorem~\ref{main theorem}(ii-iv)}$,$
and the following hold$:$

\begin {itemize}

\item[\rm(i)]
{\rm Theorem~\ref{main theorem}(v)}
holds$;$ 

\item[\rm(ii)]
{\rm Theorem~\ref{main theorem}(ii)}
takes 
$O( \xi q\log q+\mu q\log^2 \hspace{-1pt} q )$ Las Vegas
time$,$  succeeding with   probability  $>1/2$$;$

\item[\rm(iii)]
{\rm Theorem~\ref{main theorem}(iii)}
is deterministic   and takes
$O( \mu  q\log q)$  
time$,$
except in the case  $ \SU(6,q),$ where it 
takes 
$O(\xi    q\log q + \mu  q\log q)$ Las Vegas
time$,$  succeeding with   probability  $>1/2$$;$ and   

\item[\rm(iv)]
{\rm Theorem~\ref{main theorem}(iv)}
is deterministic   and takes 
$O(\mu \log q)$ time$.$

\end {itemize}

\end{theorem}

\proof 
This is  contained in  
  \cite{KS1}, except that the 
 times in   \cite[{\bf 6.6.3}]{KS1} contain a factor $q^3$
 for    the group $ \SU(6,q)$ 
 (due to the treatment of $\SU(3,q)$), which is avoided as follows.
 
It is  noted in  \cite[Sec.~5.3]{BrK2} that \cite[{\bf 4.6.3}]{KS1}
handles $\Po^-(6,q)$ in the stated times if modified using ideas  in 
\cite  {BrK2}.  This
readily gives the stated result for
$
\SU(4,q)$, which  can then be used   in
\cite  [Sec.~6]{KS1} for all larger-dimensional unitary groups.  In
particular, this leads to the stated times for $ \SU(6,q)$.  
\qed

\medskip
%\medskip

The above times do {\em not} include verification of a  presentation
of the stated groups (cf.  Theorem~\ref{main theorem}(vi));
  we will deal with that later in the context of 
of Theorem~\ref{main theorem}.
There are more recent versions of the above  theorem that run in polynomial time, assuming the availability of  suitable oracles
 \cite{BrK1,BrK2,Br2,Br3,LMO}.
 Section~\ref{Concluding remarks} contains  comments concerning  
 possible similar improvements of Theorem~\ref{main theorem}.
 We also note that  \cite{Br2,Br3} obtain
 better times than \cite{KS1} by avoiding the recursive call in the latter reference, but this has little effect on the present paper's focus on bounded rank groups. 

\convention {\rm
The proof of Theorem~\ref{main theorem}(vi) 
 uses the Las Vegas portion of Theorem~\ref{previous algorithms}(iii) for $\SU(6,q)$ when $G$ has type $\E_6(q)$, and 
the Las Vegas
 Theorem~\ref{main theorem}(vi) (more precisely, the  Theorem~\ref{main theorem}(iii) portion) for type $E_7$ when  $G$ has type $E_8(q)$; otherwise all of our uses of Theorems~\ref{previous algorithms}(iii)
 and \ref{main theorem}(iii) are deterministic.} 
 {\em In our algorithm we will avoid this {\rm``}either or{\rm''}
 possibility and assume that we are always in a deterministic 
 setting when using  the aforementioned results}. \rm
 In each   instance that is actually Las Vegas (of which there are only $O(\log q)$), up  to 20 repetitions of
 the Las Vegas version can be inserted in order to guarantee
 that  the  probability of failure is at most $1/2 ^{20}$,
which is insignificant compared to other probabilities of failure that  occur elsewhere in our algorithms.

\medskip
As in \cite{KS1,KS4,KS5,BrK1,BrK2,Br2,Br3},  we will use  crude probability  estimates, making the  number of repetitions of calls to previous routines (such as those in 
Theorem~\ref{Recognition algorithms used}(ii)) appear to be unreasonably large.  The goal has been to
prove theorems rather than to obtain best estimates  for each type of group.
 
As in \cite{KS1,KS4,KS5,BrK1,BrK2,Br2,Br3}, our algorithms contain 
statements such
as ``Choose up to $10\hspace{-.6pt}\cdot\hspace{-.6pt} 2^{12} $ pairs $z',y$\dots".  ~We could instead have used statements such
as    ``Choose~$O(1)$ pairs $z',y$\dots";
this would have eliminated some calculations, suppressed some
very  annoying constants, and looked more elegant.  However, it is not
clear how a computer would deal with such an $O(1)$ requirement.
 By contrast, ``$O(\mu \log q)$ time'' merely refers to a property of an algorithm.

% \vspace{-12pt}

\section{Groups of rank $>2$}
\label {Rank $>2$}
Throughout this section we will  assume that
\begin{eqnarray}
\label{the list} 
\lefteqn{\mbox{\em$\hat G$ is  the simply connected cover of  $F_4(q),$
$E_6(q),$ 
$\E_6(q),$  
$E_7(q)$ or 
$E_8(q)$.}}\hspace*{3.5in}
\end{eqnarray}
Here, $\hat G$ is a known copy of the group in question,
as opposed to a black box version we will eventually handle.
There are only a few cases where $\hat G$ is not also the universal cover of $G$ (cf. \cite[p.~313]{GLS}), and we will always assume that 
$q$ is large enough to avoid these.  Thus, $\hat G$
\emph{is precisely the group with that name appearing in}
 Theorem~\ref{Main Theorem}.

 We will assume the availability of the Lie algebra  
of $\hat G$.  This will be used in 
Lemma~\ref{matching} (and the Appendix),  and in Remark 
\ref{elements of hat G}.

\medskip
\noindent 
$
\begin{array}{lllllllll} 
{\bf Notation} & &&&\\
% r     &  \mbox{(twisted) rank of $\hat G$
% (for example, $r=4$ for 
%$\E_6(q)$)} \\
 \Phi  &  \mbox{root system for $\hat G$ } \\
\Phi^+ &  \mbox{the set  of positive roots } \\
\Delta &  \mbox{a base of $\Phi $} \\
p   &  \mbox{the characteristic of $\hat G$ } \\
\F  &  \mbox{$\F_q$, $q=p^e$ } \\ 
\F ' &  \mbox{$\F_{q^{\epsilon'}}$, where $\epsilon'=1$ except for $\E_6(q)$,
where  $\epsilon'=2$}
\\  
\{ f_1, \ldots, f_e \} &
 \mbox{an $\F_p$-basis  of $\F $, where $f_1=1$} \\
\{ f_1, \ldots, f_{2e} \}\ \  &
 \mbox{an $\F_p$-basis  of $\F' $ if $\hat G$ is $\E_6(q)$} 
\end{array}
$
\medskip \medskip

The {\em rank} of $\hat G$ refers to the 
twisted  rank  (for example,   
 $\E_6(q)$ has rank 4).
  
\subsection {Properties of ${\hat G}$}
\label {Properties of $G$}

We will use a standard type of   presentation for the simply connected
cover
${\hat G}$ of the simple group of Lie type we are considering.  This
presentation depends on the root system $\Phi  $ and various 
integers $C_{i,j,\a,\b},$  $ \epsilon_{\alpha  \beta },$  $ \eta_{\alpha 
\beta}, A_{\a,\b}$, all of which 
{\em we assume have been precomputed.}
%\vspace{-1pt}

\para{Presentation of the target group.}
We temporarily exclude groups of type $\E_6$. 
The  following is   just a straightforward,  shortened
version
 of the standard Curtis-Steinberg-Tits presentation \cite{St,BGKLP}.
 Use  generators $\hat X_\a(f_k)$, $\a \in \Phi$, $1 \le k \le e$, satisfying
the following relations:
\begin{eqnarray}
\label{relation1}
\lefteqn{
\mbox{\em Given any $t = \sum_k z_k f_k\in \F$  with $0 \le z_k < p$,
write
}}  \hspace*{3.5in}
 \\
\lefteqn{
\mbox{\em $\displaystyle
\hat X_\a(t) := \prod _k \hat X_\a (f_k)^{z_k};$}} \nonumber
\hspace*{2.7in}
\\
\lefteqn{\mbox{\em $
\hat X_\a (f_k)^p  =  1 $ for $ \a \in \Phi,~ 1 \le k \le e; $
}} 
\label{relation2}
\hspace*{3.5in}
\\
\lefteqn{\mbox{$\displaystyle
[  \hat X_\a (f_k), \hat X_\a (f_l) ]  =  1$ for $
         \a \in \Phi,~ 1 \le k < l\le e $; and
}}\label{relation3} \hspace*{3.5in}
\\ 
\lefteqn{\mbox{$ \displaystyle
[  \hat X_\a (f_k), \hat X_\b (f_l) ]  =
\prod_{i,j > 0} \hat X_{i\a + j\b} ( C_{i,j,\a,\b} \, f_k^i
f_l^j )~$\em  for $  \a,\b\in\Phi,~ \alpha \ne \beta,$}} 
\label{relation4} 
\hspace*{3.5in}  
\\
\lefteqn{\mbox{\raisebox{2.9ex}[10pt][0pt] {$1 \le k, l\le e.$}
}} \nonumber\hspace*{.69in}
\end{eqnarray}

\vspace{-16pt}

The right hand side of
(\ref{relation4})  is viewed as expanded, using
(\ref{relation1}), into an expression involving powers of the generators
$\hat X_\gamma(f_m)$ for $\gamma \in \Phi$, $1 \le m \le e$.
The structure constants  $C_{i,j,\a,\b}$
are  integers that are at most $2$ in absolute value (since 
we have rank $>2$), and  are given  in
\cite[Section~5.2]{Ca1}. 
The non-uniqueness of this presentation is discussed at length in \cite[p.~58]{Ca1}.
 
 %\newpage
 
The right side of  (\ref{relation4}) has at most one nontrivial term when
there is only one root length (i.\hspace{1pt}e., for types $E_6,E_7,E_8$).  In this
case, there is a nontrivial term 
$\hat X_{\a + \b} ( C_{1,2,\a,\b} \, f_kf_l )$
precisely when $\a+ \b\in \Phi$. 
A more
precise version of (\ref{relation4})  for groups of  type
$F_4$ is in   the paragraph following   (\ref{relation4'}).

The above relations provide a presentation for the simply connected cover $\hat
G$. 
\vspace{.5pt}
An algorithm for finding the center of this group is 
given in \cite[p.~198]{Ca1} using elementary linear algebra;
every element of $\zent(\hat G)$ is expressed as a word in our
generators.  However,  
$\zent(\hat G)$ can easily be found more directly 
for the groups studied here.

We will need further relations  (\ref{TG
action})--(\ref{W action on root groups}) that are consequences of the preceding ones and take  into account the action
of a split torus on the {\em root groups}
$\hat X_\alpha:=\langle\hat X_\alpha(f_k)\mid 1\le k \le e\rangle  $.

\def\E{^2\kern-.8pt E}
 \para{The group $\E_6(q)$.}
This time $\hat G  $ is the simply connected central extension of
 $\E_6(q)$,  $\Phi $ is a root system of type $F_4$ and 
$\hat G  $ has generators
$\hat X_\alpha  (f_k)$  with  $\alpha  \in \Phi $, and $1\leq k \leq e$
for
$\alpha $ long while
$1\leq  k\leq 2 e$ for $\alpha  $ short.  
We use the obvious analogues of relations
(\ref{relation1})--(\ref{relation3}), along with the  relations 
\begin{equation}
\label{relation4'}
\begin{array}{lll}
\hspace{-20pt}
[ \hat X_\alpha  (f_k), \hat X_\beta (f_l) ] = & 
\hspace{-10pt}\mbox{{\em for}:} 
 \vspace{3pt}
\\
~ 1&  \mbox{$\alpha  +\beta  \notin \Phi$} 
\\
~   \hat X_{\alpha  +\beta } ( \epsilon_{\alpha  \beta } \, f_k f_l)
&
\mbox{\em $\alpha  , \beta , \alpha  +\beta $ all short or all long} \\
~   \hat X_{\alpha  +\beta } \big ( \epsilon_{\alpha  \beta } (f_k   f_l
^q+
  f_k ^q f_l)\big  ) & \mbox{\em $\alpha  , \beta $ short$,$ $\alpha  +\beta
$ long} 
\\
~   \hat X_{\alpha  +\beta } ( \epsilon_{\alpha  \beta } \, f_k f_l)
\hat X_{\alpha  +2\beta } (
\epsilon'_{\hspace{-0.7pt}\alpha  \beta }
\, f_k f_l  f_l ^q) \ \ 
& \mbox{\em $\alpha  , \alpha  +2\beta $ 
long$,$
$\beta ,
\alpha  +\beta $ short}
\end{array}
\end{equation}
for all appropriate basis elements $f_k,f_l$.
The right hand side of
(\ref{relation4'})  is expanded as above. The structure constants
$\epsilon_{\alpha  \beta }$ and
$\epsilon_{\hspace{-0.7pt}\alpha  \beta }'$ are $\pm 1$, and as before we
assume that these have been obtained in advance. 

The relations corresponding to (\ref{relation4}) for   $F_4(q)$ are  just the relations
(\ref{relation4'})  with all  field elements in $\F$.
In this case, the third relation in (\ref{relation4'})
 involves a structure constant 
$C_{i,j,\a,\b}$  that is not 0 or
$\pm1$; and  this is the only time  this occurs for groups of rank $>2$.
 
{\em We assume that the presentations  
  $(\ref{relation1})$--$(\ref{relation4})$ or
$(\ref{relation4'})$ 
are given as part of the data describing the target group
$\hat G$.  Eventually
we will find elements of our black box group satisfying
them.}

The above presentations  are essential for  our
algorithms.  However, there are ``variants'' \cite{GKKL,GKKL2}
that may be more useful in practice:  they only involve a bounded number
of relations for any $q$ (fewer than 1000 in  \cite{GKKL}
or 50 in  \cite{GKKL2}).

\para{Additional 
relations in ${\hat G}$; the 
subgroups $T_{\hat G}$  and  $N_{\hat G}$.   }
Following \cite[p.~189]{Ca1}, if $\a \in \Phi$ and 
  $t \in \F ^*$  let
\begin{equation}
\label{Definition of $h alpha$} 
\mbox{$\hat h_{\a}(t): = \hat
n_{\a}(t)\hat n_{\a}(-1)$, \em where
$\hat n_{\a}(t) := \hat X_{\a}(t)\hat X_{-\a}(-t^{-1})\hat X_{\a}(t);$}
\end{equation}
when ${\hat G}$ is of type $\E_6$  and $\a$ is short then we also
allow $t\in \F'{} ^*$.  Define 
\begin{equation}
\label{Definition of TG}
\hspace{.2pt}
T_{\hat G}:=\langle \hat h_{\a}(t)\mid \a \in \Phi, ~t \in \F
^*\rangle  ~~{ and}~~
{N}_{\hat G}:=\langle T_{\hat G}, \hat n_{\a  }(t)\mid \alpha \in \Delta
, ~t \in
\F ^*\rangle
;\hspace{14pt}
\end{equation}
in type $\E_6$ we again use $t\in \F'{} ^*$ when $\alpha $ is short.
If 
${\hat G}$ is an  untwisted group  then $T_{\hat G}$ is a 
maximal split torus of order $(q-1)^{\mbox{\footnotesize rank of $G$}}$; if ${\hat G}$ is  
$\EH_6(q)$, then $T_{\hat G}$  is a maximally split torus of order $(q-1)^2(q^2-1)^2$.
   Moreover, $ T_{\hat G} \nor {N}_{\hat G} $, and 
   ${N}_{\hat G}/T_{\hat G}$ is the Weyl group of $\hat G $.

If $\alpha \in \Phi $ then $\hat X_\alpha $ is the set of all
$\hat X_\alpha (t)$.  
 The subgroups  $\hat X_\a$ generate  ${\hat G}$.
 
By \cite[p.~194]{Ca1},
 the root groups $\hat X_{\a}$ are 
invariant under  conjugation by $T_{\hat G}$: 
\begin{equation}
\label{TG action}
\hspace{.09pt}
\begin{array}{lll}
~\hat h_{\a}(t)\hat X_{\beta}(u)\hat h_{\a}(t)^{-1} \!=\! \hat
X_{\beta}(t^{A_{\a,\beta}}u)  &
\hspace{-3pt}
\mbox{\em except for the next instance}\hspace{-9pt}
\vspace{3pt}
\\
~\hat h_{\a}(t)\hat X_{\beta}(u)\hat h_{\a}(t)^{-1} \!=\! \hat
X_{\beta}((t t ^q) ^{A_{\a,\b}/2}  u) \ &
\hspace{-3pt}
\mbox {\em in type $\E_6,$
$\alpha\,$   short$,$
$\beta$   long$,$\hspace{9pt}%
}%     
\end{array}%
\end{equation} 
where  $A_{\a,\beta} : = 2(\a,\beta)/(\a,\a)$ for  
 the Killing  form $(~,~)$  of the underlying Lie algebra. 
 By \cite[p.~190]{Ca1} we also have
\begin{equation}
\label{W action on root groups}
\hat n_{\a}(t)\hat X_{\beta}(u)\hat n_{\a}(t)^{-1} =
 \hat X_{w_{\a}(\beta)}(\eta_{\a,\beta}t^{-A_{\a,\b}}u),
\end{equation}
where $w_{\a}$ is the reflection in the Weyl group of ${\hat G}$ corresponding to
the hyperplane  $\a^\perp $, and $\eta_{\a,\beta} = \pm 1$.
Thus, each element of  the 
Weyl group   permutes the root groups $\hat X_{\beta}$
by conjugation.

\para{How elements of $\hat G$ are described.} 
Elements of $\hat G$ are most conveniently  given in the form
$unu' $, with 
 $n\in N_{\hat G}$ and 
$u,u'$   in the Sylow $p$-subgroup  $ \langle \hat X_\gamma
(t)\mid t\in \F$  \,or\,  $ \F',\, \gamma\in \Phi ^+  \rangle$ \
(\emph{Bruhat decomposition}
 \cite[Corollary~8.4.4]{Ca1} or  \cite[Theorem~2.3.5]{GLS}).   
In this paper we do not have a natural module   as occurs in the classical   group case
\cite{CFL, KS1,Br1,Br2,Br3,BrK1,BrK2,LMO}.
However,   an element of $\hat G$ could merely be given as an automorphism of the associated Lie algebra.
See Remark~\ref{elements of hat G}  for
further discussion.

\para{Root groups and root elements.}
  The
${\hat G}$--conjugates of the $\hat X_\a$ are called {\em root groups}:  a
{\em long root group} if $\alpha $ is long and   a    {\em short root group} if
$\alpha $ is short. In case  all roots have equal length we call all root
groups ``long''.
Context will determine whether a discussion of long root groups will only be concerned with ones of the form $\hat X_\a$
rather than arbitrary conjugates of these.

Nontrivial elements of long root groups are called \emph{long root
elements}.  Each
long root element is in a uniquely determined long root group.

The following standard result is in 
\cite[Lemma~2.2]{Cooperstein}.

\begin{lemma} \label{commutator relations} 
For  long root groups $X_1,X_2,$ of  ${ \hat G},$  one of the 
following holds$:$
\begin{itemize}
\item[{\rm(i)}]
 $[X_1,X_2] = 1,$ 
\item[{\rm(ii)}]
 $|\langle X_1,X_2\rangle | = q^3$ and  $[X_1,X_2] =
\zent(\langle X_1,X_2\rangle )$ is a long root group$,$ or 
\item[{\rm(iii)}]
  $\langle X_1,X_2\rangle  \cong \SL  (2,q)$.
\end{itemize}
\end{lemma}

Two long root groups are  {\em opposite} if they generate a
subgroup isomorphic to $\SL  (2,q)$, called  a {\em long  $\SL  (2,q)$}.
Short  $\SL  (2,q)$'s are defined similarly when there are two root
lengths.
Two long root elements  are  {\em opposite} if they lie in opposite long root groups.  
Note that, when $q$ is even, two opposite long root elements 
will merely generate a dihedral group.  The preceding lemma provides a
simple way to test whether or not two long root elements are opposite: 
\begin{equation}
\label{opposite test}
\mbox{\em Long root elements $a,b$ are opposite if and only if
$\,[[a,b],a]\ne1$.\hspace{9pt}}
\end{equation}

\para{The group $\hat R  $,  the highest  root   ${\nu}$
and  the root
${\nu'}$.} Let  
\begin{equation}
\label{nu}
 \mbox{${\hat R}: = \langle {\hat X}_{\nu},{\hat X}_{-\nu}\rangle \cong
\SL(2,q),$ \em where
$\nu$ is the highest   root  of $\Phi$.\hspace{20pt}}
\end{equation}
 Then $\nu$ is a long root,   $\Delta \cup \{-\nu \}$ is the set of roots in
the extended Dynkin diagram of $\hat G$ \cite[p.~10]{GLS}, and   
\begin{equation}
\label{mu}
 \mbox{\em There is a unique long root ${\nu'}
\in \Delta$  not orthogonal to $-\nu$. \qquad\qquad\quad}
\end{equation}
Moreover, $\Delta_{\hat L}:=\Delta \cap \nu ^\perp $ is a base of the  subroot
system $\Phi_{\hat L}$ it generates, 
\vspace{1.5pt}
and $\Delta =\Delta _{\hat L}\cup
\{{\nu'} \}$.  

\para{The subgroups ${\hat L}$ and ${\hat Q}$.} 
 Define 
$${\hat L}:=\langle {\hat X}_{\a} \mid  \a \in \Phi_{\hat L}\rangle
 {~~ and }~~
  {\hat Q}:=\langle {\hat X}_{\a}\mid \alpha  \in 
\Phi^+ \setminus \Phi_{\hat L}\rangle .$$ 
If $1 \neq z \in \hat X_\nu$ then  $ \cent_{\hat G}(z) =
 {\hat Q} \semi {\hat L} =\cent_{\hat G}(\hat X_\nu)$.  The groups 
${\hat L}$ and ${\hat Q}$ are as follows:
\begin{equation}
\label{table of L}
%~~~~~~
%
\begin{array}{|l||c|c|c|c|c|l}
    \hline 
{\hat G}\raisebox{2.8ex}{~}
 & \hat F_4(q) & \hat E_6(q) & {\EH}_6(q) & \hat E_7(q) & \hat E_8(q) 
\raisebox{2.5ex} {~}
\raisebox{-1ex} {~}
\\
    \hline 
{\hat L} 
 &\Sp(6,q)& \SL  (6,q)& \SU(6,q)
 & \Spin^+_{12}(q) & \hat E_7(q) 
\raisebox{2.5ex} {~}  \raisebox{-1.3ex} {~}
\\
    \hline 
{\hat Q}
  &q^{1+14} \hspace{5pt}\flat & q^{1+20} & q^{1+20} & q^{1+32} & q^{1+56}    
\raisebox{2.5ex} {~}\raisebox{-1ex} {~}
\\
    \hline 
%\\\smallskip\smallskip
\raisebox{-2.8ex}{$\hat T_*$} \raisebox{2.8ex}{~}
  &   q^3+1  &    \frac{\displaystyle q^6-1}{\displaystyle q-1} & 
   \frac{\displaystyle q^6-1}{\displaystyle q+1}  &    q^6-1 &
\frac{\displaystyle q^{8}-1}{\displaystyle q-1}   \raisebox{3.7ex} {~}
 \raisebox{-2ex}{~}
 \vspace{-4pt}
\\
 \raisebox{2.8ex}{~}
%    \hline
  &     &    &   &   
 \hspace{-4pt} 
  \raisebox{2.8ex}{ $\frac{\displaystyle q^8-1}{\displaystyle q^{2
\raisebox{6pt}{~}\!\!}-1}$} & 
   \raisebox{2.8ex}{ $\frac{\displaystyle q^9+1}{\displaystyle q^{2
\raisebox{6pt}{~}\!\!}-q+1}$ }
\\
    \hline 
\end{array}
\end{equation}
where $ \hat E_r(q) $ denotes the 
simply connected cover of $E_r(q)$, and 
$\zent(\hat Q)=\hat X_\nu$ except where
$\flat$ indicates that this does not hold for $F_4(q)$ when 
$q$ is even (cf.  Lemma  \ref{special group lemma}(iv)).  We have also listed the orders
of  some cyclic maximal tori 
$\hat T_*$ of $\hat L$ containing $\zent(\hat G)$
that will be used  in Section~\ref{ppd
probability}.  The  orders  
in the $E_8$  case come   
from \cite[$T_{30}$ and $T_{24}$ in Table~III]{DF}.)
 
 Note that ${\hat Q}{\hat L}$ is the derived subgroup
of a parabolic subgroup $\norm_{\hat G}(\hat X_\nu)$ 
for which the unipotent radical is  ${\hat Q}$ and
the derived group of a Levi factor is   ${\hat L}$. Also, 
\begin{equation}
\label{derived group}
\mbox{\em$\zent(\hat G)<\cent_{\hat G}(\hat R) = {\hat L},$ and $T_{\hat G}$ normalizes
both
$ 
\hat R$
 and ${\hat L}$.}
\end{equation}

Define 
\begin{equation}
\label{Definition of TL}
T_{\hat L} := \langle \hat h_{\a}(t) \mid  \a \in \Phi_{\hat L},\, t \in \F^*
\rangle  \mbox{ \em \ and \ } 
{ N_{\hat L}}:=\langle T_{\hat L}, \hat n_{\a  }(t)\mid \alpha \in
\Delta_{\hat L} , ~t
\in
\F^*\rangle 
\end{equation}
with $\hat n_{\a  }(t)$ in  (\ref{Definition of $h alpha$}); 
in type $\E_6$, as in (\ref{Definition of TG})  use $t\in \F'{} ^*$ when
$\alpha $ is short, so that  
${\hat L} = \SU(6,q)$ and
$|T_{\hat L}| = (q-1)(q^2-1)^2$.
In each case,  $T_{\hat L}$ is a  maximal torus of $\hat L$,
 $T_{\hat L}< T_{\hat G}$, ${  N}_{\hat L} < { N}_{\hat G}$ and
 $ N_{\hat L}/T_{\hat L}$ is the Weyl group of $\hat L $.
 
 Let $\hat Z:=X_{\nu}$.
\begin{lemma}
\label{special group lemma} {\rm \cite[pp.~16--18]{CKS}}
\label{pairing roots} 

\begin{itemize}

\item[\rm(i)] 
 For every root $\nu \neq \a  \in \Phi^+ \setminus \Phi_{\hat L}$ 
there is a unique root $\b \in \Phi^+ \setminus \Phi_{\hat L}$ such 
that $\a + \b = \nu$.

\item[\rm(ii)] 
 If $\hat G$ is not $F_4(q)$  with   $q $ even$,$ then for each 
root group   $\hat X_\a\ne \hat Z$  in  $ \hat Q $  there is 
a unique root group $\hat X_\b $ in  $\hat Q $   that does not 
commute with $\hat X_\a$ $($and then $\a$ and $\b$ have the same length$)$.

\item[\rm(iii)] If $\hat G$ is  $F_4(q) $ with  $q $ even$,$ then for each
{
\rm long}  root group  $\hat X_\a\ne \hat Z$  in  $\hat Q $ there is 
a unique long
root group $\hat X_\b $ in  $ \hat Q $   that does not 
commute with~$\hat X_\a$.% 

\item[\rm(iv)] If $\hat G$ is  $F_4(q) $ with  $q $ even$,$ then
$\zent(\hat Q)=\langle \hat Z,\hat X_\alpha \mid \alpha $ short\/$\rangle
$  has order 
$q^7$ and is the standard module for $\hat L =\Omega (7,q)$.
\end{itemize}
\end{lemma}

 This follows from  the commutator 
relations, which also provide  more information  in the situation of this lemma:  ${\hat Q}/[\hat Q,\hat Q]$ is an $\F$-space of dimension 
$14,$  $20,$ $20,$  $32$ and $56$ in the respective cases (\ref{the list});
and it is  an  irreducible $\F {\hat L}$-module except when
$q$ is even and ${\hat G} = F_4(q)$ (producing the $\flat$ in (\ref{table of
L})), in which case
${\hat Q}/[\hat Q,\hat Q]$ has an irreducible 6-dimensional $\F {\hat L}$-submodule modulo which
it is irreducible (Section~\ref{Linear algebra in $Q/Z$} has 
computations  based on this fact).

\para{Long subgroups.}  We call any subgroup generated by (conjugates of) 
long root groups a \emph{long subgroup}.  We will especially emphasize
long subgroups such as  $\hat R$,  $\hat L$,
long $\SL(3,q)$-subgroups   such as 
$\hat S$ in (\ref{hat S})  below,  and  long subgroups
$ \Spin^ -_{8}(q)$    such as  $\hat J$ in Lemma~\ref{constructing $J$} 
below.

\para{The long subgroup ${\hat S} \cong \SL  (3,q)$.}
Let
\begin{equation}
\label{hat S}
\mbox{\hspace{-18pt} $~{\hat S} =\langle {\hat X_{ \nu}},{\hat X_{- \nu}},
{\hat X_{
{\nu'}}}, {\hat X_{-{\nu'}}} \rangle \cong \SL  (3,q),$  $T_{\hat S} = T_{\hat
G} \cap {\hat S}$ \em and
$N_{\hat S} = N_{\hat G} \cap {\hat S}$.}
\end{equation}   
The following are straightforward to check:
\begin{lemma} 
\label{lemmaS}
\begin{itemize}
\item[\rm(i)] $T_{\hat G}$ normalizes ${\hat S}$. 
\item[\rm(ii)] 
$N_{\hat G} = \langle N_{\hat S},N_{\hat L} \rangle$.
\item[\rm(iii)] If $q>2$ then 
$T_{\hat G} = \langle T_{\hat S},T_{\hat L} \rangle$
and
$N_{\hat G}/T_{\hat G}$ is the Weyl group of  $\hat G$.
\item[\rm(iv)]If $q>3$ then
$N_{\hat L} = \norm_{\hat L}(T_{\hat L}), $
$N_{\hat S} = \norm_{\hat S}(T_{\hat S})  $
and 
$N_G  =\norm_{\hat G}(T_{\hat G})$.
 
 \end{itemize}
\end{lemma}

\begin{lemma}
\label{one class of SL 3}
\label{uniquePSL}
Let $\hat S_1$ be a long subgroup of ${ \hat G}$ isomorphic to 
$\SL  (3,q)$.  Then
\begin{itemize}
\item[\rm(i)]
$\hat S_1$ is conjugate to $\hat S,$ 
\item[\rm(ii)] 
If $\hat L_1 \in L^{\hat G}$   centralizes 
  a   long
\vspace{2pt}
$\SL(2,q)$ subgroup of $\hat S_1,$
then the  pair $(\hat S_1, \hat L_1)$ is conjugate  in $\hat G$  to $(\hat S,\hat
L),$  and

\item[\rm(iii)] If ${\hat S_1}$   contains
${\hat X_{ \nu}}$   then  
$\Oh_p \big(\cent_{\hat S_1}({\hat X_{ \nu}}) \big)\le {\hat Q}$.

\end{itemize}
\end{lemma}

\proof (i) See \cite{Cooperstein} or   \cite{Liebeck-Seitz}.

\smallskip
(ii) $\hat S$ is transitive on its long  $\SL(2,q)$ subgroups.

\smallskip
(iii) 
Since the  pair $(\hat S_1, \hat X_{ \nu})$ is conjugate  in $\hat G$  to $(\hat S, \hat X_{ \nu}),$ 
we may
assume that
$\hat S_1 = {\hat S}$.  Then 
$\Oh_p \big(\cent_{\hat S}(\hat X_{ \nu}) \big)=
\hat  X_{{\nu'}} \hat X_{ \nu}\hat X_{ \nu
- {\nu'}}
\le  \hat Q$.
\qed

\begin{lemma} 
\label{L, S and tori}
$\hat G$ acts transitively by conjugation on the set of
all  $4$-tuples\break 
$(\hat L_1,  \hat S_1,  T_{\hat L_1},  T_{\hat S_1})$ with 
$\hat L_1\in \hat L^{\hat G},$ $ \cent_{\hat G}(\hat L_1)'<\hat S_1\in \hat S^{\hat G},$
 and  $T_{\hat L_1}$ and $ T_{\hat S_1}$ maximally split tori of 
  ${\hat L_1}$ and ${\hat S_1},$ respectively$,$ containing 
$\hat S_1\cap  \hat L_1$.

Moreover$,$ $\hat L$ and $\hat S\cap  \hat L$ uniquely determine $\hat S$. 
Finally$,$ if  $q>3$ then $T:=T_{\hat L_1}  T_{\hat S_1}$
is a maximally split torus of $\hat G,$ and $N/T\cong W,$ where
$N \! = \! \norm_{\hat G}(T) \! =  \! 
\langle \norm_{\hat L_1}(T_{\hat L_1}),
\norm_{\hat S_1}(T_{\hat S_1}) \rangle$.
 
\end{lemma}

\proof The preceding lemma already handles the pairs 
$(\hat L_1, \hat S_1)$.
Consider our subgroups $\hat L, \hat S$  and the 1-dimensional torus 
$\hat A:= \hat S\cap  \hat L$.  Since $\cent_{\hat G}(\hat A)$ is reductive, 
all of its maximally split tori are conjugate and contain $\hat A$.
Since $T_{\hat L}={\hat L}\cap T_{\hat G}$, this handles the triples 
$(\hat L_1, \hat S_1, T_{\hat L_1})$. 

Clearly 
 $\hat L>\cent_{\hat L}(\hat A)=\cent_{\hat G}(\hat A\hat R)
\ge \cent_{\hat G}( \hat S),$ where  
$\cent_{\hat G}( \hat S)$ is generated by 
$\zent(\hat G)$ and long root groups.
If $\hat M$ is the subgroup of $\cent_{\hat L}(\hat A)$ generated by its long root groups,  examining \cite{Ka1, Cooperstein, Liebeck-Seitz}
we find that $\hat M \zent(\hat G)=\cent_{\hat G}( \hat S)$. 
 Thus, $\hat L$ and 
$\hat A$ determine $\hat S=\cent_{\hat G}(\hat M)'$.
(In fact, $\cent_{\hat G}(\hat M)'=\cent_{\hat G}(\hat M) $ 
using \cite[Table~5.1]{LSS}.)

All maximal split tori of $\hat S$ containing $\hat A$ are $\hat R$-conjugate (as is seen by using a basis of the 3-space underlying 
$\hat S$ with respect to which $\hat A=\cent_{\hat S}(\hat R)$ consists of all matrices $\diag(\lambda, \lambda, \lambda^{-2} ))$.  Since $\hat R$ normalizes 
$\hat L, $ $\hat S$ and $ T_{\hat L}$, this proves the stated transitivity.

The final statements follow from Lemma~\ref{lemmaS}.
\qed

\para{The long subgroups  $\hat J\cong \Spin_8^-(q)$.}
\begin{lemma} 
\label{constructing $J$}
 There are long subgroups
$ \hat J \cong\Spin_8^-(q)$ containing ${\hat S}$.
\end{lemma}
\proof Each group $\hat G$ has a long subgroup $F_4(q)$  containing
$\hat S$.  Then it suffices to consider the
case  ${\hat G} = F_4(q)$,  where there is a root subsystem subgroup 
$\Spin_9(q)$ containing a conjugate of $\hat S$ that lies in a subgroup 
$ 
\Spin_8^-(q)$.~\qed

\subsection {Primitive prime divisors}
\label {Some probability estimates}
\label {primitive prime divisors}
\label {ppd probability} 
 
When the rank is $>2$,  we will always assume that $q>9$ 
in order to avoid  difficulties occurring in the  next lemma  for
small fields.   
  Remark~\ref{small q} in 
Section~\ref{Additional remarks}  discusses some of the omitted~$q$.  
 
\begin{lemma}
\label{tori probability}
\label{cyclic tori}

Let $pl$ be as  follows for the indicated types of $
\hat G$$:$
$$pl= 
\begin{cases}
p\hspace{-1pt}
   \cdot\hspace{-1pt} \ppd^\#(p;2e)\ppd^\#(p;6e)   &    F_4 \cr
p\hspace{-1pt}
   \cdot\hspace{-1pt} \ppd^\#(p;2e)\ppd^\#(p;3e)\ppd^\#(p;6e) &    E_6 \cr
p\hspace{-1pt}
   \cdot\hspace{-1pt} \ppd^\#(p;e)\ppd^\#(p;3e)\ppd^\#(p;6e) &    \E_6 \cr
p\hspace{-1pt}
   \cdot\hspace{-1pt} \ppd^\#(p;e )\ppd^\#(p;2e)\ppd^\#(p;3e)\ppd^\#(p;6e)  \  \ &
    E_7
\cr
p\hspace{-1pt}
   \cdot\hspace{-1pt} \ppd^\#(p;4e )\ppd^\#(p;8e) &
    E_7
\cr
p\hspace{-1pt}
   \cdot\hspace{-1pt}\ppd^\#(p;2e)\ppd^\#(p;4e )\ppd^\#(p;8e) &
    E_8  \cr
p\hspace{-1pt}
   \cdot\hspace{-1pt} \ppd^\#(p;2e )\ppd^\#(p;18e) &
    E_8  \cr 
    \end{cases}
$$
Let $ \ll=\ll(\hat G) $ denote the $p'$-part of $|\hat G|$. 
\begin{itemize}
\item[\rm(i)] 
 If $\tau \in \hat G$ has order of the form  $pl,$ then
$\tau^  \ll$ is a long root element or  
  $\hat G$ has type $F_4$    and $\tau^\ll$ is either a long or a short
root element. 
\item[\rm(ii)]  With probability  $\ge 1/315  q,$  an element   $\tau\in\hat G$ has
order of the form $pl$ and $\tau^  \ll$ is a long root element. 

 \end{itemize}
\end{lemma}
 
\proof
We first construct elements $\tau$ of the stated orders.
In (\ref{table of L}) we  provided  information concerning the centralizer
of both a long root element and of $\hat R$, a long root $\SL(2,q)$, 
together with   the orders of one or two maximal tori $\hat T _*$
 in that centralizer.  We will choose 
 $\tau\in \hat T _*\hat R$. 
The  integers required in the definition of $l$  exist by \cite{Zs} or the
definition of $\ppd^\#$  in Section~\ref{Background}.

These tori are constructed as follows. 
\begin{itemize}
\item  
%\vspace{-4pt}
In $\hat F_4(q) $  a subgroup $\Sp(6,q) $  centralizing a long root group  has a
cyclic maximal torus of order
$q^3+1$. 
\item 
%\vspace{-4pt}
 In  $ \hat E_6(q)$ or $\EH_6(q)$ a subgroup   
$\SL(6,q) $ or $\SU(6,q) $  centralizing a long root group  has a cyclic
maximal torus of order
$(q^6-1)/(q-1)$  or $(q^6-1)/(q+1)$, respectively.
\item 
%\vspace{-4pt}
In  $ \hat E_7(q)$  a subgroup  $\Spin^+_{12}(q) $   centralizing a long root
group   contains  subgroups $\GL(6,q)$ and $\Spin_{8}^{-}(q) \circ
\Spin_4^{-}(q)$, which produce the tori in (\ref{table of L}).   
\item 
%\vspace{-4pt}
 In $\hat E_8(q)$ a subgroup  
 $\hat E_7(q)$  centralizing a long root group  contains
subgroups of type $\SL  (8,q)\, $ (more precisely, its quotient by a central subgroup of order $(2,q-1)$)  and $\Z_{q+1}\circ {\EH_6(q)}$,
producing the   tori in (\ref{table of L}).   
\end{itemize}
\smallskip
\smallskip
(i)  By the Borel-Tits  Lemma~\cite[Theorem~3.1.3]{GLS}, 
$\tau$ lies in a 
parabolic subgroup $U\semi L$ of $\hat G$, with $ U$ unipotent containing $ \tau^ \ll  $ and $L$ a Levi factor containing
 $\zent(\hat G)$.    
Thus, we need to consider the possibility that a $p'$-element   of $L$ of order given in the lemma
 centralizes a nontrivial element of $U$.

Examination of the Levi factors that contain elements of
 order $l$
 produces the following possibilities: the normalizer of a long root group;
a parabolic of type $q^{7+8} \colon \! \hat B_3(q)$ in $\hat F_4(q)\, $ (and then $\tau ^\ll$  is a short root element); 
a parabolic of type $q^{2+6+12}\colon \!\big({\rm SL}_2(q)\circ {\rm SL}_3(q^2) \big)$ in $\EH_6(q)\, $
(and then an element of 
order  $l$
 fixes no  nontrivial  element  of the 
unipotent radical); 
a parabolic of type $q^{7+35} \colon \!\hat A_6(q)$ in $\hat E_7(q)\, $  
(and then an element of 
order  $l$
 fixes no  nontrivial  element  of the 
unipotent radical); and a parabolic of type
 $q^{8+28+56} \colon \!\hat A_7(q)$ in  $\hat E_8(q)\, $ (and then an element of 
order  $l$
 fixes no  nontrivial  element  of the 
unipotent radical). Here we used \cite{Fleischmann} in the last of these in order to verify the
statements about $\tau^\l$; references such as 
\cite{Shin, Shoj, Ca2} can also be used for other cases. 
\smallskip

   (ii)
We have $\cent_{\hat L}(\tau^p) =  \hat T_*$  in the previous
description of one type of $\tau$.  
 Also, we have
$| \norm_{\hat G}(\hat T_*\hat R)|/|\hat T_*||\hat R|\le 
 | \norm_{\hat L}(\hat T_* )\colon\!
\cent_{\hat L}(\hat T_* ) |  \le 72 $ for each of the possible tori $\hat T_*$.

Thus, there are $|\hat G\colon \! \norm_{\hat G}(\hat T_*\hat R)|
\ge |\hat G|/72|\hat T_*||\hat R|$ conjugates of
$\hat T_*\hat R$.  Even in the exceptional  $\ppd^\#$ cases (Mersenne primes, Fermat primes and $2^6-1$ in Section~\ref{Background}), each such conjugate   has at least $|\hat T_*|(1-1/2)(1-1/3)(1-1/5)(1-1/7)=
|\hat T_*|(8/35)$ elements $\tau ^p$ of the required $p'$-order 
(since $l$ has at most  four ppd-factors) and  
$|\hat R|/q$ elements of order $p$.
Thus, in  each case the number of elements  $\tau$ 
is at least $ (|\hat G|/72)({8/35})(1/q) = {| \hat G|/315 q}$.~\qed
 
 \medskip
 
 The proof shows that the probability is  better than stated.
 First of all, 2 is never a primitive prime divisor; and in all but one case there are only two or three $\ppd$-factors rather than four.  However, for simplicity we will use  the estimate $1/315 q$.

\para{Notation:}   If $ \hat G$ is of type $E_7$ or $E_8$, then there
are two choices for $l$ in the above lemma.  We will call these $l$ and  $l_0$.

\begin{lemma} 
\label{rationale for choice of tau} 
\label{generated by root groups}  
Let $\hat R_1$ be a long $\SL(2,q)$ contained in $\hat L,$
 and let $l$ $($or $l$  and $l_0)$  be as in the preceding lemma.
\begin{itemize}
\item[\rm(i)]
If $ \hat G$ is not of type $E_7$ or $E_8,$ 
%\vspace{2pt}
and if $g\in \hat L$ has order $l,$ then $ \hat L  
=\langle  \hat R_1, g    \rangle$.

\item[\rm(ii)]
If $\hat G$ is of type
$E_7$ or $E_8,$ and if  $g\in \hat L$ has order $l$ and 
$g_0\in \hat L$ has order $l_0,$  then 
$\hat L = \langle  \hat R_1,g,g_0    \rangle$. 

\end{itemize}
\end{lemma}

\proof 
Let  $\hat K:=\langle   \hat R_1^{\langle g  \rangle }   \rangle$ 
(or $\langle  \hat R_1^{\langle g,g_0 \rangle }   \rangle$ in (ii)).  
Since $\hat K$ is normalized by
$g$ (and $g_0$), as above 
the resulting ppd-factors of 
 $|\norm_{\hat G}(\hat K)| $ and the Borel-Tits  Lemma imply that 
$\Oh_p(\hat K)=1$.  Using  $|\norm_{\hat G}(\hat K)| $ 
and the lists in \cite{Ka1,Cooperstein,Liebeck-Seitz}, we see that  
$\hat K=\hat L$.~\qed
   
%%%%%%%%%%%%%%%%%%%%%%%%%%%%%%%%%%%%%%
%%%%%%%%%%%%%%%%%%%%%%%%%%%%%%%%%%%%%%

\subsection{Probability and long root elements} 
\label{Probability and long root elements} 
\label{Probability} 

Next we will study the  probabilistic behavior of
some subgroups of $\hat G$ generated by 2, 3 or 4 
 long root elements or groups.  We assume that $q>9$. 
Recall that $\hat Z =\hat X_{\nu}$.
 
%\newpage

\begin{lemma} 
\label{opposite long root element}
If $z$ is a long root element$,$ then 
 a  randomly chosen long root  element $z'$ is opposite $z$
with probability  $>1/3.$  Moreover$,$ with probability $>1/12,$
for  a  randomly chosen long root  element $z'$
either $\<z,z'\>\cong \SL(2,q)$ or $p=2$ and $\<z,z'\>$ is dihedral 
of order ${2  \ppd^\#(2 e,p)}$.
\end{lemma} 

\proof  
We may assume that $z\in \hat Z$.~The
unipotent radical
$\hat  Q = \Oh_p \big(\cent_{\hat  G}(\hat Z) \big)$ acts regularly on the set of root groups
opposite
$\hat Z$.  Then the total number of long root elements opposite $z$ is 
$(q-1)|\hat Q|$, while the total number of long root elements is $ |\hat G\colon \!
\cent_{\hat  G}(z)|$.  Hence, the desired probability is the ratio of these
integers, and it is straightforward to  check the 
lower bound  $ 1/3 $  in
all cases.

Each opposite pair $z,z'$ 
    lies in a unique long $\SL(2,q)$.  Two elements of order $p$ in that $\SL(2,q)$ generate  the required type of subgroup with probability $\ge 1/4\,$   \cite[Lemma~3.8(iii)]{KS1}. \qed

\medskip

%\newpage

  We next turn to generating the long root subgroups  
  $ \hat S=\SL(3,q$) and $ \hat J= \Spin^-_8(q)$ appearing in (\ref{hat S}) and 
Lemma~\ref{constructing $J$}.
Let $\hat R$ be as in (\ref{nu}).

Let $n( \hat S, \hat R)$ denote the number of conjugates of $\hat S$ containing $\hat R$, 
and $n( \hat J, \hat S)$ the number of conjugates of $\hat J$ containing $\hat S$.
All members of $ \hat R^{ \hat G}$ lying  in $\hat S$  are $ \hat S$-conjugate,
and all members of $ \hat S^{ \hat G}$ lying in $ \hat J$  are $ \hat J$-conjugate.
Therefore, the numbers $n({\bf X},{\bf Y})$, $ ({\bf X},{\bf Y})=( \hat S, \hat R) $  or $( \hat J, \hat S)$,
can be  obtained from
Tables~\ref{n(S,R)} and \ref{n(J,S)}  by 
simplifying 
the obvious formula to 
$$n({\bf X},{\bf Y})=\frac{| \norm_{\hat G}({\bf Y}) |{\hspace{.4pt}}
\raisebox{-.9ex} {\hspace{.00001pt}} |{\bf X} |}
{  | \norm_{ \hat G}({\bf X}) |{\hspace{.4pt}} | \norm_{\bf X}({\bf Y}) |  }.
$$

\begin{table} 
\caption{Number of root $\SL(3,q)$  that contain  
 a given long root 
 $\SL(2,q)$}

\vskip 0.5cm
\label{n(S,R)}
\begin{center}
\begin{tabular}{|c|c|c|}
\hline 
$\hat G$ & $n(\hat S,\hat R)$ & $q$-exponent 
\raisebox{3ex} {~} \raisebox{-.7ex} {~}
\\
\hline
\hline
$G_2(q)$ & $\displaystyle\frac{q(q+1)}{2}$ 
\raisebox{4ex} {~} \raisebox{-2.6ex} {~}
& 2 \\
\hline
$\D_4(q)$ & $\displaystyle\frac{q^3(q^3+1)}{2}$
\raisebox{4ex} {~} \raisebox{-2.6ex} {~}
 & 6 \\
\hline
$F_4(q)$ & $\displaystyle\frac{q^6(q^3+1)(q^4-1)}{2(q-1)}$ 
\raisebox{4ex} {~} \raisebox{-2.6ex} {~}
& 12 \\
\hline
$E_6(q)$ & $\displaystyle\frac{q^9(q^3+1)(q^2+1)(q^5-1)}{2(q-1)}$
\raisebox{4ex} {~} \raisebox{-2.6ex} {~}
 & 18 \\
\hline
$\E_6(q)$ & $\displaystyle\frac{q^9(q+1)(q^3+1)(q^5+1)}{2}$ 
\raisebox{4ex} {~} \raisebox{-2.6ex} {~}
& 18 \\
\hline
$E_7(q)$ & $\displaystyle\frac{q^{15}(q^3+1)(q^5+1)(q^8-1)}{2(q-1)}$ 
\raisebox{4ex} {~} \raisebox{-2.6ex} {~}
& 30 \\
\hline
$E_8(q)$ & $\displaystyle\frac{q^{27}(q^9+1)(q^5+1)(q^{14}-1)}{2(q-1)}$
\raisebox{4ex} {~} \raisebox{-2.6ex} {~}
 & 54 \\
\hline 
\end{tabular}
\end{center}
\end{table}

\begin{table}  
\caption{Number of root $\Spin_8^-(q)$  that contain a given long root $\SL(3,q)$ }

\vskip 0.5cm
\label{n(J,S)}
\begin{center}
\begin{tabular}{|c|c|c|}
\hline
$\hat G$ & $n(\hat J,\hat S)$ 
\raisebox{3ex} {~} \raisebox{-.7ex} {~}
& $q$-exponent \\
\hline
\hline
$F_4(q)$ & $\displaystyle\frac{q^3(q^3-1)}{2}$ 
\raisebox{4ex} {~} \raisebox{-2.6ex} {~}
& 6 \\
\hline
$E_6(q)$ &$\displaystyle\frac{q^6(q^3-1)^2}{2}$ 
\raisebox{4ex} {~} \raisebox{-2.6ex} {~}
& 12 \\
\hline
$ \E_6(q)$ &$\displaystyle\frac{q^6(q^3+1)^2}
{\raisebox{1.4ex} {~}2}$ 
\raisebox{4ex} {~} \raisebox{-2.6ex} {~}
& 12 \\
\hline
$E_7(q)$ &$\displaystyle\frac{q^{12}(q^6-1)(q^5-1)(q^3-1)}
{\raisebox{1.6ex} {~}2(q^2-1)}$ 
\raisebox{4ex} {~} \raisebox{-2.6ex} {~}
& 24 \\
\hline
$E_8(q)$ &$\displaystyle\frac{q^{24}(q^{12}-1)(q^9-1) (q^5-1)}
{\raisebox{1.6ex} {~}2 (q^2-1)}$ 
\raisebox{4ex} {~} \raisebox{-2.6ex} {~}
& 48 \\
\hline
\end{tabular}
\end{center}
\end{table}

\begin{lemma} 
\label{S J probabilities}

Let  $\hat R,$ $\hat S$ and $\hat J$ be as before. 
\begin{itemize}
\item [\rm(i)]  
The probability  is 
at least $1/3$ that $\hat R,$ together with a conjugate 
of $\hat Z$ opposite $\hat Z,$ generate  a conjugate of $\hat S$.

\item [\rm(ii)] 
The probability  is   at least $1/3$  that $\hat S,$ together with a conjugate 
of $\hat Z ,$ generate  a conjugate of $\hat J$. 
 \end{itemize}
\end{lemma}

\proof  
For $ ({\bf X},{\bf Y})\hspace{-1.6pt}=\hspace{-1.6pt}
(\hat S,\hat R) $  or $(\hat J,\hat S)$, 
the desired probability is  at least 
 $   n({\bf X},{\bf Y})\b /|\hat Q|$, where  $\b $ is the 
number of conjugates $\hat Z'$ of $\hat Z$
 inside {\bf X} that are  opposite $\hat Z$ and satisfy  ${\bf X} = \langle {\bf Y},\hat Z'\rangle$
(recall that $|\hat Q|$ is the number of $\hat G$-conjugates 
of $\hat Z$   opposite $\hat Z$). 
From Tables~\ref{n(S,R)} and \ref{n(J,S)} we obtain 
${n(\hat S,\hat R)}/{|\hat Q|} \geq  {4q^{-3}}/{9}$ and 
${n(\hat J,\hat S)}/{|\hat Q|} \geq  
q^{-9}(1-q^{-3})^2/2.$
It remains to estimate $\b $ in our two cases. 
\smallskip

(i) 
Let $V$ be the natural module for $\hat S=\SL(3,q)$. Then 
$V = [V,\hat R] \oplus \cent _V(\hat R)$, and 
the only maximal overgroups of $\hat R$ in $\hat S$ are 
the parabolics 
$\norm _{\hat S}([V,\hat R])$ and 
$\norm _{\hat S}\big(\cent _V(\hat R) \big)$. 
If $\hat Z'<\hat S$ is a 
conjugate of $\hat Z$ opposite $\hat Z$, then 
$[V,\hat R]  \ne \cent _V(\hat Z')$  
and $\cent _V(\hat R)\ne [V,\hat Z']$.  

Thus, if also  $\hat S>\< \hat R,\hat  Z'\>$, then either 
$[V,\hat R] > [  V, \hat  Z'] $ or 
$  \cent _V(\hat R)  < \cent _V(\hat Z') $. 
There are at most $2q^2$ such $\hat Z'$ out of the $q^3$ opposite $\hat Z$. Thus $\b \geq q^3 - 2q^2 $,
and the desired probability is at least  $({4q^{-3}}/{9} )(q^3 - 2q^2 ) >1/3$.
\smallskip

(ii) Let $V$ be the natural 8-dimensional module for 
$\Omega^-(8,q)$; we will view all subgroups of $\hat J$ as subgroups of 
$\Omega^-(8,q)$.  Then $\hat S$ splits $V$ as 
$V={V_6^+\perp V_2^-}$.  Long root groups $\hat Z'$ correspond to totally singular 2-spaces $T $ of $V$ via  $T=[V, \hat Z']$.  If $V= \langle T, V_6^+ \rangle $  then 
$T^\perp\cap V_2^-=0$ (as otherwise $T$ and $V_6^+$ would lie in a 7-space).  Consequently,  
$\langle \hat S,\hat Z' \rangle $ is an irreducible subgroup 
of $\hat J$ generated by long root groups and hence is $\hat J$  
\cite{{Ka1,Liebeck-Seitz}}.  

Thus, we only need to estimate the number of totally singular 2-spaces {\em not} spanning $V$ together with $V_6^+$.  Each such 2-space contains a point of $V_6^+$.  Since $V_6^+$ has $(q^2+q+1)(q^2+1)$ singular points, and each is contained in $(q^3+1)(q+1)$ totally singular 2-spaces, 
there are at most $(q^2+q+1)(q^2+1)(q^3+1)(q+1)<3q^8$ totally singular 2-spaces meeting $V_6^+$ (as $q > 9)$. 
There are   $q^9$ long root groups in $\hat J$ opposite $\hat Z$. It follows that 
$\b \ge  q^9-3q^8 ,   $
so that the desired probability is at least  
$(q^9-3q^8 )\cdot q^{-9}(1-q^{-3})^2/2\ge 1/3$.~\qed

\medskip
We will need   variations on the previous   arguments:

\begin{lemma} 
\label{S J probabilities2}
\begin{itemize}
\item [\rm(i)]  
   Suppose that   $D$ is a subgroup generated by 
opposite long root elements $z,z'  $  
such that either $D\cong \SL(2,q)$ or $q$ is even and $D$ is dihedral of order $2 \ppd^\#(2e,2)$.
Then the probability 
 is at least $1/4$
that $D,$ together with a 
conjugate $y$ of $z$  
 opposite $z,$ generate a conjugate of $\hat S$.
 \item [\rm(ii)]  
 The
probability is at least $1/3$ that 
$\cent_{\hat S}(\hat R)= \hat S\cap \hat L$  and a random conjugate $\hat
S^l,$ $l\in\hat  L,$ generate  
a $\hat G$-conjugate of $\hat J$ having an element normalizing $\hat R$ and conjugating $\cent_{\hat S}(\hat R)$ into $\hat S^l$.

 \end{itemize}
 \end{lemma}

\proof 
(i)   By Lemma~\ref{S J probabilities}(i),  
with probability  at least  ${1}/{3}$ the three root groups containing $z,z'$ or $ y$
generate a conjugate of  $\hat S$. Thus,  we only need  a lower bound on the 
conditional probability that   $\hat S =\<D,y\>$ for 
a root element $y\in \hat S$
opposite $z$.

 In view of the structure of  $D$, the 
only maximal overgroups of $D$ in $\hat S$ are    $\norm_{\hat S}([V, D ])$ and 
$\norm_{\hat S}(\cent_V(D)  )$ (compare \cite[Lemma~3.7]{KS1}). 
Define $\b$ as  at the start of the proof of Lemma~\ref{S J probabilities}. Then at least  $\b (q-1)$ of the $q^3(q-1)$ root elements  in $\hat S$ opposite   $z$ generate $\hat S$ together with $D$, so that  the desired probability is at least 
$(1/3)\b (q-1)/q^3(q-1)\ge (q^3 - 2q^2)/3q^3> 1/4$. 

%\newpage
 
 \smallskip
 
(ii) Somewhat as in Lemma~\ref{S J probabilities}(ii),
\begin{equation}
\label{$S^l$}
\mbox{\hspace{-16pt}\em$\hat S$ and $\hat S^l$ generate a long  
$\Spin^- _8(q)$ subgroup
with  probability  $ > 1/4 $.  }
\end{equation}
For, the number of ``good'' conjugates 
$\hat S^l$ such that $\<\hat S, \hat S^l\>\in\hat J^{\hat G}$   
is $n(\hat J,\hat S) 
\hspace{-.15pt}\cdot\hspace{-.15pt}\g$,~where
$n(\hat J,\hat S)$ is in Table~\ref{n(J,S)} and 
 $\g$ is  the number of  
 good $\hat S^l$ per $\hat J$-conjugate containing $\hat S$. 
 On the other hand, $|\hat S^{\hat L}|$ is just the number 
$n(\hat S, \hat R ) $
 of conjugates 
  of $\hat S$ containing $\hat R$ (by Lemma~\ref{one class of SL 3}). Thus, the 
 desired probability is 
 $n(\hat J,\hat S)\g/n(\hat S,\hat R)$.
We will provide a lower bound for $\g$, from which (\ref{$S^l$}) will follow using Tables~\ref{n(S,R)} and \ref{n(J,S)}.

For this purpose, we restrict our attention to the $8$--space $V$ associated with  $\hat J$.
As for the preceding lemma, $\hat S^l$ is good if (and only if) 
$\<\hat S  ,\hat S^l\>$ is irreducible on  $V$.

If $V_6^+: =[V, \hat S]$  then $V_6^+=U_3\oplus U_3^*$ for 
totally singular $\hat S$-invariant 3-spaces $U_3, U_3^*$. 
Also, $ [V, \hat R]$  is a
nondegenerate 4-space of type $+$ meeting $U_3, U_3^*$
at 2-spaces $U_2, U_2^*$, respectively.

Since $l$ centralizes $\hat R$, the totally singular 3-spaces $U_3^l, U_3^*{}^l$ meet 
$ [V, \hat R]$ in totally singular 2-spaces 
 $\tilde U_2:=U_2^l, \tilde U_2^*:=U_2^*{}^l$,
lying in the same ``half'' of the set of totally singular 2-spaces of 
  $ [V, \hat R]$  as $U_2, U_2^*$ (this ``half'' consists of $q+1$ totally singular 2-spaces pairwise meeting in 0,
  all of which are $\hat R$-invariant);  there are $(q+1)q $ such ordered pairs $\tilde U_2, \tilde U_2^* $ of
distinct 2-spaces. 
Each of the subspaces $U_3^l, U_3^*{}^l$    meets the $4^-$-space 
$[V, \hat R]^\perp=\<\tilde U_2, \tilde U_2^*  \>^\perp$ in a singular point; there are $(q^2+1)q^2$
  ordered pairs $p_1,p_2$   of distinct singular points in $[V, \hat R]^\perp$.
Each such ordered pairs  of
 2-spaces and of points determine   a unique  ordered  pair
$  \< \tilde U_2 , p_1 \>,$ $   \<  \tilde U_2^* , p_2 \>$ of totally singular 3-spaces left invariant by a conjugate $\hat S^l$, and each $\hat S^l$ arises this way exactly twice 
(twice because the ordered pairs $\tilde U_2 , \tilde U_2^* $
and $p_1,p_2$ determine the same conjugate of $\hat S$ as the ordered pairs 
 $\tilde U_2^*, \tilde U_2$
and $p_2,p_1$).  Of the pairs $p_1,p_2$ of singular points,  
$q^2-1 $ points $p_1$ do not lie in $V_6^+$ and at least 
$q^2-1 -(q-1) $ points $p_2$ do not lie in $\<V_6^+,p_1\>$,
 in which case $V=\<V_6^+,p_1,p_2\>=\<U_3, U_3^*,U_3^l, U_3^*{}^l\>$ and 
$\<\hat S, \hat S^l\>$ is irreducible.  Thus, 
$2\g\ge (q+1)q\hspace{-1pt}
   \cdot\hspace{-1pt} (q^2-1)(q^2-q )$.
Now Tables~\ref{n(S,R)} and \ref{n(J,S)}  yield  (\ref{$S^l$}).
\smallskip

Let $\hat A:=\hat S  \cap  \hat L$.
{\em It remains to show that $\langle \hat S^l, \hat A\rangle\cong\hat J,$
$l\in L,$  assuming that
$\langle \hat S, \hat S^l\rangle \cong\hat J$.}  
Instead of this it will be more convenient to show that
{\em $\langle \hat S, \hat  A^l\rangle =\hat J,$
$l\in L,$  assuming that
$\langle \hat S, \hat S^l\rangle =\hat J$.}
 
If $\langle \hat S, \hat A^l\rangle$ is irreducible on $V$ then   
so is $\langle \hat S^{\langle \hat S, \hat A^l\rangle} \rangle$, and then  both of these groups are $\hat J$ using \cite{Ka1,  Liebeck-Seitz}. 
Moreover, $\hat S$ and $\hat S^l$ are long $\SL(3,q)$-subgroups of
$\hat J$ containing $\hat R$, 
and hence are conjugate under $\norm_{\hat J}(\hat R)\:$ (cf. Lemma~\ref{one class of SL 3}(ii)).  Then $\hat A^l{}^j<\hat S^l{}^j= \hat S$ for some $j\in   \norm_{\hat J}(\hat R)$, as required in (ii).

We will   assume  that 
$\langle \hat S, \hat A^l\rangle$
is reducible and obtain a contradiction.
A generator  %$c$ 
of $\hat A^l=\cent_{\hat S^l}(\hat R)$ acts on
$V$  
 by  centralizing a $2^-$--space (hence with eigenvalue 1 there) while acting on $V_6^+{}^l$ with two invariant totally singular 3-spaces
  $U_3^l, U_3^*{}^l$ and eigenvalues on them of the form
$\lambda, \lambda, \lambda ^{-2}$  and  
$ \lambda ^{-1}, \lambda ^{-1}, \lambda ^{2}$, respectively, 
where $\lambda$ has order $q-1$.
In particular, $\cent_{V_6^+{}^l}(\hat A^l)=0$ since $q>3$,
so that $\cent_{V}(\hat A ^l)$ has no singular points.
Thus,   the only 
 totally singular 3-spaces left invariant by $\hat A^l$ are 
contained in $V_6^+{}^l$. 
%,

Any proper $\langle \hat S,\hat  A^l\rangle$-invariant subspace 
$W$ of smallest dimension must be totally singular or nondegenerate. 
Clearly $\hat S$ and $\hat A^l$ have no fixed common nonzero vector
since 
$\cent_V(\hat S^l)=\cent _V(\hat A^l)$ and 
$\hat J=\langle \hat S, \hat S^l\rangle$.  Thus, $W$ is $U_3$ or $ U_3^*$, and yet we have seen that it must be   contained in $V_6^+{}^l$.
Then the 6-spaces $ V_6^+ =[V, \hat S]$ and $V_6^+{}^l=[V, \hat A^l] $ both contain  both
$W $  and the $4^+$-space $[V, \hat R]$, and those span at least a 5-space. Thus,
$\<[V, \hat S],[V, \hat S^l]\> = \<[V, \hat S],[V, \hat A^l]\> <V$ and 
$\langle \hat S, \hat S^l\rangle$ is reducible.  This is the desired contradiction.~\qed

 %\newpage

\subsection {Start of the proof of Theorem~\ref{Main Theorem}}
\label {Start of proof}

We are given a black box group $G$ that is a nontrivial epimorphic image of the
universal cover  $\hat G$ of an  exceptional group of Lie type of
rank $ >2$ over a field of  order $q>9$.
Therefore $\hat G$ is the simply connected cover  
\cite[p.~313]{GLS}.  We start by using  the Monte
Carlo algorithm in
\cite{BKPS} in order to (probably) find the type of group we are dealing
with.  Similarly, every time we call an existing constructive
  recognition algorithm
  in Theorem~\ref{previous algorithms}
   we assume that \cite{BKPS} has first been 
used in order to make it likely that we are testing a group having the  desired structure: 
the algorithm in 
\cite{BKPS} is far faster than any constructive recognition algorithm (such as Theorem~\ref{previous algorithms}),
although these checks are not necessary for the proof of
Theorem~\ref{main theorem}.

Eventually we will test that the group is, indeed, as expected:  
in Proposition~\ref {constructing lambda}
and Corollary~\ref {verification}
we
will verify  a presentation 
(\ref{relation1})--(\ref{relation4}) or
(\ref{relation4'})   for $G$. Such a
presentation   is also crucial for uses of  Theorem~\ref{main theorem},
such as those   in
\cite{KS2, LG}.

\subsection {Finding a  long  root element}
\label {Finding a root element}  

Choose up to $  3150q$ elements $\tau\in G$ until one is found such that 
$|\tau| =p l$ for $l$   in  Lemma~\ref{tori probability}.
 When we obtain  $\tau$ of  the desired sort,  Lemma~\ref{tori probability}(i)  states that $z:=\tau^\ll$ is a long root element,
or possibly
a short one when $G$ has  type  
$F_4$.  For the latter groups we proceed somewhat differently.

Suppose that $G$ has type $F_4$.
If $q$ is even then the graph 
automorphism sends short root elements to long ones, so  we may assume
that $z $ is long.  If $q$ is odd then we  run the algorithm
 up to $3200q$  times,
from finding $\tau$ until the group $L$
is constructed and tested at the start of Section~\ref{Finding  $L$}  
(specifically: we find $\tau $ and then find and test  
$z',y ,
S, S_2, Z_1, J$ and $L$).  

\Remark Correctness\rm: There is
no subgroup of $F_4(q)$, $q$ odd, generated by short root elements  and
isomorphic to $\Sp(6,q)$.  For,  $F_4(q)$ has exactly 2 classes of involutions, with centralizers
$\Spin_9(q)$ and   $\big(\SL(2,q)\circ\Sp(6,q) \big)\hspace{-1pt}
   \cdot\hspace{-1pt}2$ (for a long $\SL(2,q)$)
    \cite{Shoj};
only the latter type has a subgroup $\Sp(6,q)$, and 
the  long
root groups in $\Sp(6,q)$ are also long for  $G$. 
Thus, if we obtain a subgroup $L\cong \Sp(6,q)$ then $z$ is a long root element.

There are other ways to handle this odd case.  For example,
the group generated by 4 conjugates of a long root element is
 isomorphic to $\Spin_8^-(q)$ with probability $\ge 1/16$, but the same is not true for short root elements, once again using the nature of the centralizers noted above of the 2 classes of involutions.
  In Section~\ref{Concluding remarks}, Remark~\ref{Involution
centralizers.},  this ambiguity is avoided   using an entirely different
approach that finds the involution in $R$ and then its centralizer in $G$.
\medskip

For the cases $E_7(q)$
and $E_8(q)$ there are two possibilities    $l,l_0$
in Lemma~\ref{tori probability}, and  hence we also
find a second element $\tau_0 $  of   order $pl_0$.
Then $z_0:=\tau^\ll$ is a long root element.

\rel $\ge 1-1/2^9$ for $\tau$  and $\tau_0 $ in all but the
exceptional  $F_4$ case. For, 
all $\tau$ fail to have the required  order with probability $\le(1-1/315q)^{3150q}<1/2^{10}$, by Lemma~\ref{tori probability}.

In the exceptional $F_4$  case,  for a given choice $\tau$, if $z$ is a long root element then 
we will succeed at showing this and finding the needed
elements and  subgroups with probability $\ge 1-1/2^{8}$ (in view of the individual probabilities in the next sections).
Hence, we will succeed for a given $\tau$ with probability 
$\ge (1/315q)(1-1/2^{8})> 1/320q$.  All $3200q $
repetitions fail with probability $<(1-1/320q  )^{3200q}<1/2^{10}$.

\timing
$O(q[\xi + \mu   \log ^2\hspace{-1pt}q])$ to choose elements $\tau$ (and  $\tau_0$) and to  test  the
order of each of them using \cite[Lemma~2.7]{KS1}; but 
$O(\xi q\log q + \mu q\log^2\hspace{-1pt}q)$  if the $F_4$ test is needed.

\subsection{Matching up root elements.}
For the cases $E_7(q)$ and $E_8(q)$
we have two elements $ \tau $ and $\tau_0 $, and we have powers $z$ and
$z_0$ of them that are  long root elements.  
{\em We need to arrange to have $\<z\>=\<z_0\>$.}
 
Repeat up to $2 40$ times: choose a conjugate $z_1$ of $z_0$,
test whether $z$  and $z_1$  are opposite; and for odd $q$ 
use Theorem~\ref{previous algorithms}(ii) to test whether
$\langle
z,z _1\rangle\cong\SL  (2,q)$,   and, if so,  to
obtain a 
constructive isomorphism $ \SL(2,q)\to \langle z,z_1\rangle$.
If  $p=2$ then
$\langle z,z_1\rangle $ is dihedral of order dividing $2(q\pm1)$.

For each   $q$  it is now easy to conjugate 
 $\<z \>$ to $\<z_1\>$ and hence to $\<z_0\>$. 
  
Thus, {\em we can conjugate $\tau_0$ in order to 
arrange    that}  
$\<z\>=\<\tau ^\ll\>= \<\tau_0^\ll\>=\<z_0\>$
(recall that $ \ll $ denotes the $p'$-part of $|\hat G|$).

\rel  $\ge1-1/2^{10}$. For, 
by Lemma~\ref{opposite long root element}, $z_1$ is opposite $z$   with probability
$>1/3$, and we use   (\ref{opposite test}) to  test this. 
If this
occurs  and  $q$ is even then we are merely conjugating within a dihedral group.

If $q$ is odd then  
$\langle z,z_1\rangle\cong  \SL(2,q)$ with probability $\ge 1/12$ 
(by Lemma~\ref {opposite long root element}), in which event   Theorem~\ref{previous algorithms}  succeeds with probability $>1/2$. 
 Thus, all $2 40 $ repetitions 
 fail with probability $<(1-1/24)^{2 40}<1/2^{10}$.

\timing $O(\xi q\log q + \mu q\log^2\hspace{-1pt}q),$ dominated by the time to find the constructive isomorphism.  

%\newpage 
 
\subsection {The subgroups $R$, $Z$, $Z^-$   and $S$}
\label {Finding a root group}
\label {The subgroups $R$, $Z$, $Z^-$   and $S$}
 
 Choose up to $10\hspace{-1pt}\cdot\hspace{-1pt}  2^{12}$ pairs $z', y $ of conjugates     of $z$, and  
use (\ref{opposite test})  and
Theorem~\ref{previous algorithms}(ii)
in order to test whether both are  opposite   $z$ 
and  $S:=\langle z,z',y\rangle$  and 
$ S_2:=\langle z, z' {}^{\tau ^p },y
\rangle$ are both isomorphic to ${\hat S}=\SL  (3,q) $; and, if so, 
\vspace{.5pt}
    to find constructive isomorphisms
$\Psi _S\colon    \!{\hat S} \to S   $  and 
$\Psi_{S_2} \hspace{-1pt}\colon \!{\hat S} \to S_2$, together with generating sets $\SSS_{\hat S}$  and $\SSS_S^*$  of $\hat S $ and $S   $, respectively, such that $  \SSS_{\hat S}\Psi_S=\SSS_S^*$.  
{\em We may assume that ${\hat S}$  is the subgroup of ${\hat G}$  defined in 
\rm(\ref{hat S}); \em we will use the notation in} (\ref{nu}).

Find $R:=\hat R \Psi_{S } <   S$,
$Z:=\hat X_\nu \Psi_S$ and $Z^-:=\hat X_ {-\nu}\Psi_S$
using Theorem~\ref{previous algorithms}(iv).
 Then  $R =\langle
Z,Z^{-}\rangle \cong  \SL (2,q)$.

Conjugate within  ${\hat S}$  in order to have   
$z\in Z $ and $z'\in Z^- $. Then $\tau^p$ centralizes  $Z$ since it centralizes $z\in Z$.
Find the  root group   $Y<S$ containing  $y$.

Use $\Psi_{S_2}$   to find  an element  of 
 $\Oh_p\big (\cent_{S_2}(Z)\big )$ conjugating    $ (Z^{-}) {}^{\tau ^p}$  to $Y$
 (recall that  $Z$  and  $Y$ are opposite), and
use $\Psi_{S }$ to  find an element  of
$\Oh_p\big (\cent_{S}(Z)\big )$ conjugating  $Y$  to $ Z^{-} $.
The  product of these two elements is  an element
$c\in \Oh_p \big(\cent_{s_2}(Z) \big)  \Oh_p\big (\cent_S(Z)\big ) \subseteq
Q:=\Oh_p\big (\cent_G(Z)\big )$    such that 
 $(Z^{-} ){}^{\tau ^pc} =Z^{-} $ (cf.~Lemma \ref{one class of SL 3}(iii); of course we do not yet have $Q$ to work with).  Then
$\tau ^pc$ and  $\tau ^p$ are elements 
of $\cent_G(Z)$ that agree mod~$Q$,
so that $l$ divides the order of
$\tau ^pc$.  Moreover, $\tau ^pc$ normalizes $Z^-$
while centralizing $Z$.

Thus, {\em $\tau^pc$ centralizes $  R$ and has order divisible by $l$.}

Recall from the preceding section  that $Z$ contains $\<z\>=\<z_0\> $ when $G$ is of type $E_7$ or
$E_8$.  In that case  we have a second element $\tau_0$,
and  we obtain in the same way a
second element  $\tau_0^pc_0$  of $\cent_G(R)$, this time 
{\em of  order divisible by $l_0$.} 
 
 \rel  $\ge1-1/2^{9}$. 
 For, both members of a pair $z', y$ are opposite $z$, with  $z'$ behaving as in the second part of Lemma~\ref{opposite long root element} and
$S, S_2\cong \SL(3,q)$, with probability 
 $>(1/12)(1/3)(1/4)^2 >1/2^{10} $  (by Lemmas~\ref{opposite long root element} and  \ref{S J probabilities2}(i)); in which case   
 Theorem~\ref{previous algorithms}(ii) succeeds for both $S$ and $S_2$  with probability $> (1/2)^2$.  
Hence,   all
  $10\hspace{-.2pt} \cdot\hspace{-.2pt}  2^{12} $  
 repetitions  fail with
probability $\,<(1-1/ 2^{12})^{10\cdot 2^{12} }<1/2^{10}. $ The probability   involved in 
repeating the above for $\tau_0$, if needed, is dealt with similarly.

\timing 
$O(\xi q \log q + \mu q\log ^{2}\hspace{-1pt}q)$, dominated by 
finding   
$\Psi _S$ and $\Psi_{S_2}$
using Theorem \ref{previous algorithms}(ii). 

 \subsection{The long subgroups  $J$ and    $R_1$.} 
\label{The long subgroups  $J$ and    $R_1$.} 
Repeat up to $30$  times:
choose a  conjugate $Z_1$ of $Z$, and use   Theorem~\ref{previous algorithms}(ii)   in order to
test whether $J:=\langle S, Z_1\rangle\cong \Spin^{-}_8(q) $; and, if so, to obtain a constructive isomorphism $\Psi_J\colon \! \Spin^{-}_8(q) \to J$.%
  
Find a long $\SL(2,q)$-subgroup   
$R_1<\cent_J(R)$ using  $\Psi_J $.  {\em
Obtaining this long $\SL(2,q)$  is the only use we have  for
$J$ and  $\Psi_J$}.

\rel  $\ge1-1/2^{10}$ using Lemma~\ref{S J probabilities}(ii).  
 
\timing 
$O(\xi q \log q + \mu q\log ^{2}\hspace{-1pt}q)$, dominated by 
finding  $\Psi
_J$ using Theorem~\ref{previous algorithms}(ii). 

\subsection{The subgroups $L$, $T$ and $N$}
\label{Finding  $L$}

Let
$L:=\langle  R_1,  \tau^pc   \rangle  $ 
or $\langle  R_1,\tau^pc,\tau_0^pc_{{}_0}\rangle  $
in the cases 
$F_4,$ $ E_6,$ $ ^2\hspace{-1pt}E_6$ or
$E_7, $ $E_8$, respectively.  The
generators of $L$ lie in 
in $\cent_G(R)$
(cf. Section~\ref{Finding a root group}).
Hence,   
${L=\cent_G(R)}$  by  Lemma~\ref{rationale for choice of tau}.

  The subgroups
$S$ and $L$ behave  as in Lemma~\ref{one class of SL 3}(ii), and hence the pair
$S, L$ is uniquely determined up to conjugacy in $G$.  In particular,  we can
use the information in Section~\ref {Properties of $G$} to study $G$   by means
of  constructive  isomorphisms for  these subgroups.  Note, however, that
these isomorphisms might not  match up properly, which will make us
(possibly) have to modify the pair $(S,L)$ in Lemma~\ref{matching}.

Use  
up to $10$ repetitions of   
Theorem~\ref{previous algorithms}(ii), 
or recursion if $G=E_8(q)$,
in order  to find  
generating sets 
$\SSS^*_L$ of
$L$ and
 $\hat \SSS_L$ of $\hat L$ and   an isomorphism
 $\Psi _L\colon \!
\hat L\to L$  such that 
$\hat \SSS_L\Psi _L=\SSS^*_L$.
Also find  the
following  subgroups of $G$
using  (\ref{Definition of TL}) and (\ref{hat S}): 
$$
  T_L : =  T_{\hat L}\Psi_L ,~
 T_S: = T_{\hat S}\Psi_S ,~
 N_{ L}: =N_{\hat L}\Psi _L \mbox{~\em and~} 
 N_{ S} : =N_{\hat S}\Psi _S. 
$$
(Recall that we already have  a generating set $\SSS^*_S$ of $ S$.)
We will often use the fact that $\Psi_S$  and $ \Psi_L$ are isomorphisms even though the target epimorphism 
$\Psi=\Psi_G$ may not be bijective. In particular,  $ \Psi_L^{-1}$ always produces a unique element of $\hat G$.

\rel $\ge 1-1/2^{10}$.

\timing 
$O( \xi q \log q+ \mu  q \log ^{2}\hspace{-1pt}q)$,  dominated by
finding $\Psi_L$. 

  \Remark \rm
\label{matching presentations}
A version of the presentation (\ref{relation1})--(\ref{relation4}) or
(\ref{relation4'}) is used   for $L$ 
as part of  Theorem~\ref{previous algorithms}(ii).
 Conceivably this is not a subpresentation of the presentation 
 (\ref{relation1})--(\ref{relation4}) or
(\ref{relation4'}) that we are using for $\hat G$:
the signs may not agree.
 We {\em assume} that,
as part of the  recursive call, the signs in the presentation 
 (\ref{relation1})--(\ref{relation4}) or
(\ref{relation4'}) for $\hat L$ have been changed so as to coincide with the corresponding 
ones for $\hat G$.    Since we are only dealing with presentations of 
groups of small (bounded) rank,  
there are only a few sign changes required here.
%\medskip

\subsection{Matching up $T_S$ and $T_L$ in order to obtain $T$.} 
\label{Matching}
At this point it need not be the case that  $\<T_S, T_L\>$ is a maximal torus of $G$.  In order to guarantee that property
we need to arrange for the  1--dimensional torus 
$S\cap L$  of
both
$S$ and $L$ to lie in both of the tori  $T_S$ and $T_L$:
 
 \begin{lemma}
\label{matching}
There is an algorithm  replacing the pair  $(S,L)$ by a conjugate  pair  in order to have  
$ S\cap L  = T_S\cap T_L$.  This algorithm is deterministic and runs in     $O(\mu q)$ time$,$ 
except when $G$  is $E_8(q),$  in which case it is Las Vegas$,$
takes $O(\xi q \log q + \mu q   \log^2\hspace{-1pt}q  )$
 time and succeeds  
with probability  $\ge 1-1/2^{10}$.

 \end{lemma}

 \proof  
 Recall from Section~\ref{Finding a root group} that  ${\hat S}$  is the subgroup of ${\hat G}$  defined in  (\ref{hat S}). 
 Since  $R=\hat R \Psi_{S }$, we can 
find $S\cap L = \cent_S(R)= \big(\cent_{\hat S} (\hat R)\big)\Psi_S$  using 
 Theorem \ref{previous algorithms}(iv). 
 Since $T_{\hat S}$ normalizes the root groups 
 ${\hat X}_{\nu},{\hat X}_{-\nu}$  of $\hat R$,  
$ T_{\hat S} $  contains 
$ \cent_{\hat S} (\hat R)  $
(using a basis of the 3-space for $\hat S$ as in the proof of Lemma~\ref{L, S and tori}).
Thus,  $S\cap L =\big(\cent_{\hat S} (\hat R)\big)\Psi_S \le 
  T_{\hat S}\Psi_S=T_S .$

We will provide two entirely different approaches to the remaining part of the proof:  arranging to have  $S\cap L\le
T_L$. The first is deterministic (as in the statement of the lemma) and simpler for $G$ not of type $E_8$,
while the second is more uniform.
The timing in the lemma refers to the first method.
(For rank 2 groups in Section~\ref{Rank $2$ groups} 
we will use the first method.)%
  \smallskip

{\noindent \bf Method 1.}
 We assume initially  that  $G$ does not have type $E_8$. Then $\hat L$ is (essentially) a classical group 
(cf. (\ref{table of L})); let $V$ be its    natural module.
(It will not matter that   this module is not faithful when
$\hat L$ is a spin group.) 

We have found the (cyclic) group  $S\cap L$ using $S$.  
Find $\hat A:={(S\cap L)\Psi_L^{-1}}$ 
using   Theorem~\ref{previous algorithms}(iii). 
 Diagonalize $\hat A$ on 
$V$ using  a    hyperbolic basis  that  determines a maximal  split torus $\hat T$ of $\hat L$ containing $\hat A$.  
Find $\hat l$ in the classical group $\hat L$  
such that $\hat T^{\hat l}=T_{\hat L}$
(this is just  a basis change).
Find $l:=\hat l\Psi_L$ using  Theorem~\ref{previous algorithms}(iv).
Replace $S$ by $S^l$ and $T_S$ by $T_S^l$.
\ %
({\em Correctness}:   We have   $S^l\cap L=(S\cap L)^l= \hat A^{\hat l}\Psi_L< \hat T^{\hat l}\Psi_L=T_{\hat L}\Psi_L $, where the latter  is $ T_L$ by definition in Section~\ref{Finding  $L$}.
Then  $S^l\cap L= (S\cap L)^l\le T_S^l\cap T_L \le S^l\cap L $
since $S\cap L\le T_S$.  Therefore, replacing $S$ by $S^l$ and $T_S$ by $T_S^l$  gives the desired equality 
$S\cap L = T_S\cap T_L$.)

If $G$  has type $E_8$ we again  find
$\hat A:= ( S\cap  L)\Psi_L^{-1} $, using up to 10 recursive calls to Theorem~\ref{Main Theorem}(iii,vii).
Then  the following are accomplished in the Appendix:
find the Lie algebra of   $\hat L  \cong \hat E_7(q)$,  
 then find a Chevalley basis    
producing a split torus of $\hat L$ containing $\hat A$, and finally find  $\hat l\in \hat L$ conjugating this torus to the torus 
$T_{\hat L}$ in (\ref{Definition of TL}).  
Find $l:=\hat l\Psi_L$ using another recursive call to Theorem~\ref{Main Theorem}, and replace $S$ by $S^l$ and $T_S$ by $T_S^l$.
\ %
({\em Correctness}:  Once again $S^l\cap L \le  T_{\hat L}\Psi_L =  T_L   $, and our replacement again gives
$S\cap L =   T_S \cap T_L .$)

 % \newpage 

   \smallskip
{\noindent \bf Method 2.}
Find  the subgroup 
$A:=(\hat S \cap \hat L)\Psi_L$  of  $ \hat T_{\hat L}\Psi_L=T_L$
using Theorem~\ref{previous algorithms}(iii).

Repeat  up to $30 $ times: choose $l\in L$,   use
Theorem~\ref{previous algorithms}(a) to test whether
$\langle S ^l, A\rangle \cong \Spin_8^-(q)$ 
 and, if so, use the resulting  constructive
isomorphism 
$\Spin_8^-(q) \to \langle S^l , A \rangle $
in order to find $j\in \langle S^l , A \rangle$   that
normalizes $R$ and conjugates $A $ into   $  S^l$. 
Let $m:= lj  ^{-1}$.
Replace  $S$ by  $S^{m}$ and
$T_S$ by  $T_S^{m}\!$.

{\em Correctness}: 
There is an epimorphism $\Psi\colon\! \hat G\to G$ 
extending $\Psi_L$ and hence sending $\hat R$ to $R$. 
 Then $\hat S\Psi$  contains
$\hat R\Psi=R$, and $A=(\hat S \cap \hat L)\Psi_L
=(\hat S \cap \hat L)\Psi = \hat S \Psi\cap \hat L \Psi 
= \cent_{\hat S \Psi } (R)  $
behaves as in Lemma~\ref{S J probabilities2}(ii).

 By that lemma, 
we may assume that $\langle S ^l, A\rangle$ is isomorphic to 
$\Spin_8^-(q)$ and has an element  
normalizing $R$ and 
conjugating $A$ into $S^l$.
With $m\in \norm_G(R)$ as above, $A\le S^{m}\cap L$,
so that  $A= {S^{m}\cap L }$ 
by Lemma~\ref{one class of SL 3}(ii)
since $ A=\hat S \Psi\cap   L  $.
Then  $A =(S\cap L)^{m}   <T_S^{m}  $  
  (by the start of the proof of this lemma),  while
  $A<T_L $   by definition.  Thus, 
  $A\le T_S^m\cap T_L\le S^m\cap L= A$.
 Replacing $S$ by $S^{m}$   and
$T_S$ by  $T_S^{m}$ gives $T_S\cap T_L= S\cap L $.

\timing  
Method~2 requires
$O(\xi q \log q + \mu q   \log^2\hspace{-1pt}q  )$ time,
dominated by the test for isomorphism with $\Spin_8^-(q)$.
 
 Method~1  uses  
 Theorem~\ref{previous algorithms}(iii,iv)  for $\Psi_L$,
 and hence runs in   $O(\mu q\log q)$ time  if  $G$ does not have type
$E_8$.  However, in the $E_8$ case it again  takes 
$O(\xi q \log q + \mu q   \log^2\hspace{-1pt}q  )$ time
since a constructive isomorphism is used in the Appendix.
  (N.\hspace{-.1pt}B.--The faster $O(\mu q\log q)$ time is significant,
  but it does not influence the overall time for the algorithm in
Theorem~\ref{main theorem}.)

\rel $\ge 1-1/2^{10}$ in Method 2, in view of 
Lemma~\ref{S J probabilities2}(ii) and the 
   $30$ repetitions of Theorem~\ref{previous algorithms}(ii). 
The same probability can be obtained    in the $E_8$ case
of  Method 1.
\qed
 %\newpage

\smallskip
\smallskip

At this point we could also arrange to have 
$\Psi_S|_{\hat  S\cap \hat L} = \Psi_L|_{\hat 
S\cap \hat L}  $,  but  we will not need this.  

%  \bigskip 
 
\para{The subgroups $T$, $N$ and  $W$.}\:
 By  Lemmas~\ref{L, S and tori} and \ref{matching}, $T:=\<T_S,T_L\>$ is a maximal torus    and    
 $W:=N/T$ is the Weyl group of $G$, where  $N : = \langle N_{ S},N_{ L}\rangle $.  
 
 \subsection{The root groups $X_\alpha$}
\label{The root groups $Xalpha$}

Associated with $W$ there is a root system $\Phi $   having a subsystem
$\Phi _L$ corresponding to  $L$.
In Section~\ref{The subgroups $R$, $Z$, $Z^-$   and $S$}  we already used the roots $\nu, \nu'$ appearing in
(\ref{nu}, \ref{mu})
since 
$\hat S$ was defined using (\ref{hat S}).
There is a base $\Delta _L$  for $\Phi _L$ 
such that 
$ \Delta:= \Delta _L \cup    \{{\nu'} \}$ is  a base for $\Phi  $.  

We next find the $|\Phi|$ root groups $X_\a$, $\a\in \Phi$.  
We already have  $X_{ \nu}=Z $  and $X_{-\nu}=Z^{-} $.   
Use  the isomorphism  $\Psi_L$  and 
Theorem~\ref{previous algorithms}(iv) 
%(or a recursive call in the $E_8$ case) 
to find the $T_L$-invariant root groups $X_\alpha, \alpha\in \Phi_L $.
 Conjugate these   using $N$ 
 in order to obtain all $|\Phi|\le 240$ root groups $X_\alpha, \alpha\in \Phi$.%

\timing $O(\mu \log q )$ using  $\Psi_L$
(Theorem~\ref{previous algorithms}(iv)).
%, or
%$O(\xi q \log q + \mu q   \log^2\hspace{-1pt}q  )$
%for the recursive call in the $E_8$ case.  
For, we only need one nontrivial root element in one root group $X_\a$ of each length,  an element $h_\a(t)$ 
generating the corresponding 1-dimensional torus, and
a ``reflection'' $ n_\b(1)$ for each $\b\in \Delta_L$, after which we
can conjugate using    (\ref{TG action})  and (\ref{W action on root
groups}).

  Note also that we only need coset representatives in $N$ of the
stabilizer  in $N $ of the long root 
$\nu$;  this  stabilizer   is  $N_LT $.
A similar remark holds for  short roots, if there are any.  
 There are at most 240 such coset representatives for each type of root, and these
can be quickly found in $O(1) $ time using standard permutation group algorithms for $W$
\cite[Ch.~4]{Ser}.  Alternatively, it is straightforward to 
write  coset  representatives as explicit products of fundamental reflections in  the Weyl group.
 
\subsection{The epimorphism  $\Psi\colon \! \hat G \to G_0$}
\label{Finding and verifying the presentation}
Let  $G_0:=\langle  X_\a\mid \a\in \Phi\rangle$.
We next  show that $G_0$   is an epimorphic image of $\hat G$.
In  Corollary~\ref{verification} we will test whether  
each member of the  original generating set $\SSS$ of $G$ 
lies in
$G_0$, thereby verifying that $G_0$   is $G$.

The isomorphism $\Psi_L$ 
lets us ``coordinatize'' each root group $X_\alpha, \alpha\in  
\Phi_L$:  labeling the elements of $X_\alpha$  as   
$ X_\alpha(t) 
,  t\in \F$ or $\F'$, in a manner
preserved by  the   conjugations (\ref{W action on root groups}) for   
$\a \in \Phi_L$ and satisfying the relations
(\ref{relation1})--(\ref{relation4}) or (\ref{relation4'}).
This was already noted in  Remark~\ref{matching presentations}.
 We need to   coordinatize  each root group
$X_\alpha,\alpha\in  \Phi$, in the same manner:

\begin{proposition}
\label{constructing lambda}\ 
There is a deterministic $O(\mu  \log^2\hspace{-1pt} q)$-time    algorithm that  labels any given element of any root group
$X_\a, $  $\a\in \Phi,$  as $X_\a (t)$ for some $t$ in
$\F$ or  $\F',$ in such a way 
that the map $ \hat X_\a (f_k)\mapsto X_\a (f_k)$ $($for all appropriate $\alpha $  and $k)$
extends to an epimorphism $\Psi\colon \! \hat G\to G_0 $.
Moreover, $\Psi | _{\hat L}=\Psi_L$.

\end{proposition}

\proof
We have $\hat G$ and its presentation, and we  have 
already coordinatized all 
$X_\alpha(f_k)=\hat X_\alpha(f_k) \Psi_L, \alpha\in \Phi_L$.

We also already have the long root $\nu'$  in  (\ref{mu}).  
By (\ref{TG action}), 
  $\norm_{\hat R}(\hat X_ {\nu'}) $ centralizes $\hat L$ and  is transitive
on the nontrivial elements of
$\hat X_ {\nu'} $.
Hence, we can choose any nontrivial element of $X_ {\nu'}  $  and label it
  $X_{{ {\nu'}}}(1)$.  
  We now show that \emph{all remaining labels are uniquely determined.}

 Let $\delta\in \Delta_L$  be the long root not perpendicular to ${\nu'}$. 
Using (\ref{TG action})  for $h_\delta(f_k)$
we can correctly label 
$X_{{ {\nu'}}}(f_k)$ and hence any given element of $X_{{ {\nu'}}}$.

By $\!$(\ref{relation4}),~we have relations 
$[\hat X_\alpha(f_k ),\hspace{-.9pt} \hat X_\beta(f_l) ]\! =\! \hat
X_{\alpha +\beta}(\epsilon_{\alpha, \beta } f_k f_l)$ 
in  $\hat G$   whenever $\alpha\in \Phi_L$,
$\beta $ and $\alpha + \beta$   are long.
Each subgroup  $X_\a$ of $L$ has already been coordinatized.   
 Starting with all root groups of $L$ together
  with $X_{\b}:=X_{\nu'}$,
by repeatedly using these relations with hats removed 
we coordinatize all positive  long root groups. 
Alternatively, we could achieve  this by using  (\ref{W action on
root groups}) for $n_\b(1), \b\in \Phi_L$.

We next coordinatize   $X_{-{\nu'}}$   using
$\a=\nu-{\nu'} \in \Phi_L$, 
$\beta= - {\nu'}$ together with the desired relation
$[ X_{{\nu'}+\alpha}(1),  X_{-{\nu'}}(u) ]
 =  X_{\alpha}(\epsilon_{
{\nu'}+\alpha, -{\nu'}\: }u )$ in (\ref{relation4}) or
(\ref{relation4'})  (here $\epsilon_{{\nu'}+\alpha, -{\nu'} }
 := C_{1,1,{\nu'}+\alpha, -{\nu'}}$ in (\ref{relation4})).
First, find an $\F_p$-basis for the elementary abelian group $X_{-{\nu'}}$ 
(recall that this root group was obtained as a conjugate of a root group of $L$).  For each element $x$  in this basis, 
find its coordinate $u$ via $[ X_{{\nu'}+\alpha}(1),  x ]
 =  X_{\alpha}(\epsilon_{{\nu'}+\alpha, -{\nu'} \:}u )$ using linear algebra in $X_\a$.
This produces the coordinates of a basis of $X_{-{\nu'}}$
and hence of any given element of $X_{-{\nu'}}$. 

Now coordinatize all negative long root groups as above.
 
This leaves us with the groups of type $F_4$ or $\E_6$, where there 
are also short roots to consider. Here we use the last relation in 
(\ref{relation4'}) as above in order to coordinatize  $X_{\alpha + \beta}$ whenever 
$\alpha  , \alpha  +2\beta $  are long and 
$\beta \in \Phi_L,\alpha  +\beta $ are short.
Namely
$[ X_\alpha  (1), X_\beta (f_l) ] =
   X_{\alpha  +\beta } ( \epsilon_{\alpha  \beta } \,  f_l)
X_{\alpha  +2\beta } 
(\epsilon'_{\hspace{-0.7pt}\alpha  \beta }
 f_l^{\hspace{0.7pt}}  f_l ^q)
$
where we already know $X_{\alpha  +2\beta } 
(\epsilon'_{\hspace{-0.7pt}\alpha  \beta }
 f_l^{\hspace{0.7pt}}  f_l ^q)$.

Finally, we   verify all of  the relations 
  $(\ref{relation1})$--$(\ref{relation4})$ or
$(\ref{relation4'})$.
This proves our assertions concerning both $\Psi$ and $\Psi|_{\hat L}$.

This algorithm is deterministic.  The stated time includes verifying the  relations  (cf.   \cite[{\bf7.2.2}]{KS1}). \qed 
  \medskip

Note that this same commutator method could have been used to  {\em produce}
all of the root groups $X_{\alpha}$, not just to label them.  This may, in
fact, be more efficient in practice.
Also note that $\Psi$ extends $\Psi_L$  but not necessarily   $\Psi_S$.  

   \Remark  \rm
 \label{S*}   
We have $G_0  = \langle \SSS^*\rangle ,$ 
\vspace{1.15pt}
 where  $\SSS ^*$ consists of  all of the  $ X_\alpha (f_k),$
 $\alpha\in\Phi$.
Let $\hat \SSS$ consist of  the elements $\hat X_\alpha (f_k)$
of $\hat G$, so that $\hat \SSS\Psi \! =\! \SSS^*$ is the defining property of $\Psi$.%

\begin{corollary}
\label{random G0}
A random element of $G_0$ can be constructed as a 
straight-line program of  
length $O(\log q)$ in
$\hat \SSS$   in time $O(\mu\log q)$.

\end{corollary}

\proof
Let $\displaystyle \hat U:=\prod_{\a>0} \hat X_\a$ and 
$\displaystyle \hat U_w:=\! \!\prod_{\a>0>w(\a )} \! \!\hat X_\a$ for each 
$w\in W=N_{\hat G}/T_{\hat G}$ (for a suitable order of the factors).
Also let $h_\d$ be a generator of $\hat h_\d(\F^*)$
(or of $\hat h_\d(\F'{}^*)$ if $\d$ is short), for 
  each $\d\in \Delta$.
Then $T_{\hat G}$ is the direct product of the groups $\<h_\d\>, \d\in \Delta.$
For $w\in W$ choose $n_w\in N_{\hat G}$ such that 
$w=n_wT_{\hat G}$.

By \cite[Corollary~8.4.4]{Ca1} or  \cite[Theorem~2.3.5]{GLS}, 
every  element of $\hat G$ has the unique 
Bruhat normal form 
$unv$ with 
$u\in \hat U, n\in N_{\hat G}, 
w:=n T_{\hat G}\in  W$ and $ v\in \hat U_w$.

Hence, a random element of $\hat G$ is obtained by choosing $w$ and  hence $n_w$,  then $t\in T_{\hat G}$ and hence 
$n:=n_wt$, and finally letting $u$ and $v$ be products of randomly chosen elements of the relevant root groups.
By (\ref{relation1}), each of the  $O(1)$ root group elements 
appearing in the definition of $\hat U$ or $\hat U_w$ is a product of powers of elements of $\hat \SSS$ with exponents between $0$ and $p-1$, hence 
can be obtained using a straight-line program
of length $O(\log q$) from $\hat \SSS$. Similarly,
$t=\prod_{\d \in \Delta}h_\d^{a(\d)}$ with $0\le a(\d)<|h_\d|$,
and  (\ref{Definition of $h alpha$})
shows that $t$ also
can be obtained using a straight-line program
of length $O(\log q$) from $\hat \SSS$.
 Thus,  the required 
 random  root group  elements and $t$
 are obtained by randomly choosing $w$ and all of the preceding exponents. 
 
Finally, apply $\Psi$   in order to obtain
 a random element of $\hat G\Psi=G_0$.~\qed
 
 \medskip
 
 Note that this corollary involves the more classical notion of ``random" element rather than the more subtle version in \cite{Bab} (cf. Section~\ref{Background}).
 In particular, the parameter $\xi$ is not involved.

\subsection {Effective transitivity of $Q$}

The set $Z^G$  of long root groups is far too large to be managed effectively 
using standard permutation group methods   (cf. \cite{Ser}).
Nevertheless, as in \cite{KS1,Br2,BrK1,BrK2,LMO}, we need to circumvent this
difficulty when using   the action of  $Q:=\langle 
X_\alpha\mid 
\alpha\in 
\Phi^+
\setminus
\Phi_L\rangle$ on this set. 
     As in the above references, the following 
\emph{effective transitivity} of $Q$ will be crucial later
 (in  Section~\ref
{Straight-line programs}):

 \begin{lemma}
\label{transitivity of $Q$} \label {Transitivity of $Q$}
{\hspace{-1pt}}
There is an $O(\xi q \log q+ \mu
q\log ^2\hspace{-1pt}q)$--time 
 Las Vegas algorithm which$,$
with probability $>1-1/2^{10},$
 when given long root groups $A$ and $B$ opposite
to $Z,$  finds the unique element $u\in Q$ such that $A^u=B$.

\end{lemma}

\proof 
 Each  long root group opposite $Z$ has the form $B^v$ for a unique
 $v\in {Q}$.
Repeat up to  $60$ times:  choose $v \in Q$,    test whether 
$S(v):= \langle Z,A,B^v\rangle\cong\SL  (3,q)$ using   
 Theorem~\ref{previous algorithms}(ii); if so 
obtain  a constructive isomorphism $\Psi_{S(v)}\colon \!
 \SL(3,q)\to
Y$, and finally use $\Psi_{S(v)}$ and Theorem~\ref{previous algorithms}(iii,iv) in order to obtain an element of $\Oh_p\big(\cent_{S(v)}(Z)\big)$ conjugating $A$ to $B^v$. 
Since $ \Oh_p \big (\cent_{S(v)}(Z)\big)$ is transitive on $A^Q\cap {S(v)}$,
such an element exists, and it is in $Q$ by 
 Lemma~\ref{one class of SL 3}(iii).

\rel $\ge 1-1/2^{10}$, since 
$ \langle Z,A,B^v\rangle\cong\SL  (3,q)$
with probability $\ge 1/3$ by 
Lemma~\ref{S J probabilities}(i),
and Theorem~\ref{previous algorithms}(ii) succeeds  with probability
 $>1/2$,  so that all $60$ repetitions fail with probability 
$\le (1-1/6)^{60}<1/2^{10}$.

\timing $O(\xi q\log q + \mu q\log^2\hspace{-1pt}q)$, dominated by
   finding $  \Psi_{S(v)}$. 
\qed

\subsection {Linear algebra in $Q/Z$}
\label {Linear algebra in $Q/Z$} 

We  next address the problem of 
writing an element $g \in Q$ as a word 
in the generators $X_\a(t)$.   
 
  Fix an ordering of the roots for $Q$, with $Z=X_
\nu$ first.  (For example, modify  the ordering in \cite[p.78]{Ca1} so
that $
\nu$ is first.)   Then each  $g \in Q$ can be written as a product $g =
\prod _{\a\in\Phi^+\setminus
\Phi_L} X_\a(t_\a)$ in the chosen order,
with each $t_\a\in \F$ or $\F'$ written as $\F_p$-linear combinations of the given bases of $\F$ or $\F'$. 
We will  call this product    the  {\em  standard form} of $g$. 

\begin{proposition}
\label{expressing}
\label{Expressing}
The standard form of any given   $g \in Q$    can be computed 
deterministically in $O (\mu\log q)$ time. 
\end{proposition}

\proof 
We first deal with the case  in which
{\em   $G$ is not   $F_4(q) $ with $q$ even.}
(The omitted case is handled in  the following lemma.) 
We must find the standard form
$\prod _{\a \in \Phi^+ \setminus \Phi_L} X_\a(t_\a)$
  of $g$. Let $X_{\gamma}(t_{\gamma})$ be the 
rightmost nontrivial factor in the product. By Lemma~\ref{pairing
roots}(ii) there is a unique root group $X_\b$ in $Q$ that does not 
commute with   $X_{\gamma}$. 
Then we can find $t_\gamma$ using linear algebra in $X_{\nu}$:
$$
[g,X_{\b}(1)] = [X_{\gamma}(t_{\gamma}),X_\b(1)] = 
X_{\nu}(C_{\gamma,\b,1,1}t_{\gamma})
$$ 
by (\ref{relation4}) and (\ref{relation4'}),
since  $X_{\beta}$ commutes with $g_1 := g X_{\gamma}(-t_{\gamma})$.

Now compute $g_1$ and repeat $O(1)$ times. The procedure   ends with $g\in X_\nu=Z$ after  
we have processed  $O(1) $ roots in $\Phi^+ \setminus \Phi_L$.  

This procedure is deterministic.   The time takes into account  the need to write a given field element  $C_{\gamma,\b,1,1}t_{\gamma}$ in terms of the basis vectors $f_k$.

\medskip

\para{The case $F_4(q)$ with $q$ even.}
  Here we will  modify the above
procedure using explicit knowledge of the positive roots of the root system of type $F_4$ together with the explicit presentation 
(\ref{relation1})--(\ref{relation4}) or
(\ref{relation4'}). 

{\bf Conventions:}  The roots in our base $\Delta $ are 
ordered $\alpha _1,\alpha
_2,\alpha _3,\alpha _4$, so that the high root  is $\nu
=2342$, where   we   write
$\alpha =abcd$ if 
$\a = a\a_1+b\a_2+c\a_3+d \a_4$.   

  {\bf The positive roots:}

  $ {1000}, {0100},{0010},{0001},{1100},{0110},{0011}
,{1110},{0120},{0111},{1120},{1111} ,$

$ {0121}
,{1220},{1121},{0122}
,{1221},{1122},  {1231},{1222}, {1232},{1242},{1342},{2342}$

{\bf The roots for $Q$:}   those   of the form 
$1bcd$ or 2342.

{\bf The short roots for $L$:}  $\pm0001,~ \pm0011,~ \pm0010, ~\pm0110,~
\pm0111,~
 \pm0121$.

{\bf The long roots for $L$:} $\pm0122,~\pm 0120 , ~ \pm0100$.

{\bf The short roots for $Q$:} 1232, 1231, 1221, 1121, 1111,  1110. 

{\bf The long roots $\ne 2342$ for  $Q$:}  1341, 1242, 1222, 1122, 1220, 1120,
1100, 1000. 

The above lists of $n=6$ or 8 roots in $Q$ are listed so that    the $i$th and
$(n-i+1)$st roots sum to the highest root.  For example, 1231 + 1111 =
2342 and 1222 + 1120 = 2342.   

\begin{lemma}
\label{peel F4}
{\rm Proposition~\ref{expressing}} holds if $G$ is $F_4(q)$
with $q$  even. 
\end{lemma}

\proof We must find the standard form  
$\prod _{\a \in \Phi^+ \setminus \Phi_L} X_\a(t_\a)$ of $g$. As all short
root groups of 
$Q$ lie in the center of $Q$ we can move all long root factors of $g$ to the
end (the  right hand side) of the product, and then   compute the long 
root ``coordinates" as above for the root groups $\ne X_\nu$. 

It remains to find the standard  form of an element   
$g \in \zent(Q)= 
\langle  X_{2342}, X_{1232},$ $ X_{1231}, X_{1221}, X_{1121},
X_{1111},  X_{1110}     \rangle $.  We repeatedly  use (\ref{relation4})
for these short root groups.

Compute $s_0:= [[g,X_\a(1)],X_{-\a_4}(1)]$, where $\a = 0121 $.
By (\ref{relation4}),
$s_0 = X_{1232}( t_{1110})$, from which we find  $t_{1110}$.

Define $s_1:=gX_{1110}(t_{1110})$ and  compute   
$[s_1,X_{0111}(t)] = X_{1232}( t_{1121})$  in order to find $t_{1121}$.

Define $s_2 :=  s_1X_{1121}(t_{1121})$ and  compute   
$[s_2,X_{0121}(t)] = X_{1232}( t_{1111})$ in order to find $t_{1111}$.

Define $s_3 :=  s_2X_{1111}(t_{1111})$ and compute   
$[s_3,X_{0011}(t)] = X_{1232}( t_{1221})$in order to find  $t_{1221}$.

Define $s_4 :=  s_3X_{1221}(t_{1221})$ and compute   
$[s_4,X_{0001}(t)] = X_{1232}( t_{1231})$ in order to find  $t_{1231}$.

Define $s_5 :=  s_4X_{1231}(t_{1231})$ and compute 
$[s_5,X_{-0001}(t)] = X_{1231}( t_{1232})$ in order to find $t_{1232}$.

Finally, compute  $ s_5X_{1232}(t_{1232})=X_{2342}(t_{2342})\in Z=X_{2342}$  in
order to find $t_{2342}$.

Once again this procedure is deterministic and the time is clear.
\qed

%\vspace{-4pt}

%\newpage

\subsection{Straight-line programs; testing that $G=G_0$}
\label{Straight-line programs}
We can now prove parts (ii) and (iii)  of Theorem~\ref{main
theorem}.
First of all we may need to slightly increase the set 
  $\SSS^*$ in Remark~\ref{S*} in order to use recursion. 
In Section~\ref {Finding $L$}  we used either Theorem~\ref{previous algorithms}(ii), or a recursive call when $G$ is
$E_8(q)$,  in order to find a new generating set  $\SSS^*_L$ for $L$.
If necessary, increase  $\SSS^*$ 
by adjoining this set, in which case  adjoin 
  $\hat \SSS _L=\SSS^* _L\Psi_L^{-1}$  to $\hat \SSS$
  (cf. Remark~\ref{S*}).
Thus,  $\hat  \SSS$ and $ \SSS^*$ still have size
$O(\log q)$ and $\hat \SSS\Psi = \SSS^*$.
This takes  $O(\mu q \log q)$ time by Theorem~\ref{previous
algorithms}(iii).

 \begin{proposition}
\label{SL P}
\begin{itemize}
\item[\rm(i)] \hspace{-5pt}
\vspace{1pt}
There is a deterministic $O(\mu \log q)$--time algorithm  which$,$ when
\vspace{1pt}
given   $\hat g\in \hat G,$
finds  $\hat g\Psi $ and a straight-line program of length $O(\log
q)$ from
$\hat \SSS $ to   $\hat g$.   

\item[\rm(ii)]  
There is a deterministic $O(\mu \log q)$--time  algorithm that finds 
a generator of $\zent(\hat G)$. 
\item[\rm(iii)] 
There is an $O(\xi q \log q+\mu q \log ^ 2\hspace{-1pt}q  )$--time  
Las Vegas algorithm which$,$ with probability $\ge 1-1/2^7,$
when  given $g\in G$  finds a preimage $g\Psi^{-1}$  and a
straight-line program of length $O( \log q )$ from $ \SSS^*$ to $g$. 

\end{itemize}

\end{proposition}

\proof
(i)
We have assumed the availability of the Lie algebra for $\hat G$ and the action of  $\hat G$ on that algebra.
Use \cite[Theorem~8.1]{CMT} and \cite{CHM} to write $\hat g$ in the Bruhat  form $unu' $, with  $n\in N_{\hat G}$ and  
$u,u'$   in the Sylow $p$-subgroup  $ \langle \hat X_\gamma
(f_k)\mid $  all appropriate $k$ and $\gamma\in \Phi^+ \rangle$. 
Then use  (\ref{relation1})--(\ref{relation4}) or
(\ref{relation4'}), together with  
(\ref{TG action})--(\ref{W action on root groups}), in order
to write   $u,u'$ and $n $
in terms of straight-line programs from $\hat \SSS\,$   
  \cite{Riebeek,CMT,CHM} (compare Corollary~\ref{random G0}). 
  Apply $\Psi $ in order to obtain a straight-line program from $\hat \SSS\Psi $ to   $\hat g\Psi $.
 
\smallskip 
(ii) 
There is an algorithm in \cite[pp.~198--199]{Ca1} for finding 
$\zent(\hat G)$.  However, for  each of  the present small number of exceptional groups
(\ref{the list})  one can instead readily write down the center of
$\hat G$ in terms of the elements $\hat h_{\alpha _i}(t)$,
and hence in terms of the elements  $\hat X_{\pm \alpha _i}(f_k)$, in 
$O(\log q)$ time (cf. (\ref{Definition of $h alpha$})).
Now the center of $G$ is obtained using (i). 
  \smallskip 

(iii)  
Use Corollary~\ref{random G0}
to choose up to 30   elements ${y\in G_0}  $  
in order
to find one such that $[[z^{gy},z],z]\ne1$,
so that $Z$ and $Z^{gy}$ are opposite  by (\ref{opposite test}).

Find $u\in Q$ such that
$Z^{gyu}=Z^-$ using Lemma~\ref{transitivity of $Q$};  find  a
straight-line program  of length $O(\log q)$ from $\SSS^*$ to $u$ using  Lemma~\ref{Expressing}. 
Now $gyun_\nu  $ normalizes $Z$, where $n_\nu: =n_\nu (1) $ 
  is defined  using (\ref{Definition of
$h alpha$}) without the hats.  
It follows that the desired result holds for $g$ if it holds for  $gyun_\nu 
$. 

 Thus, 
we will replace $g$ by $gyun_\nu  $, so that $g$ normalizes $Z$.  Now
$Z^-{}^g$ is opposite $Z$.  Again use 
Lemmas~\ref{transitivity of $Q$} and  \ref{Expressing}
in order to find $u'\in Q$ such that $Z^-{}^g{}^{u'}=Z^-$, as well as 
 a straight-line program of length $O(\log q)$ from $\SSS^*$ to $u'$.
Thus, we may now assume that $g$ normalizes both $Z$ and $Z^-$.

Find $h=h_{\nu'}(t)$ acting on
$Z$ and $Z^-$ in the same manner as $g$,
using    (\ref{mu}) and  (\ref{TG action}).
Find a straight-line program of length $O(\log q)$ from 
$\SSS^*$ to  $h^{-1}$ using (\ref{Definition of $h alpha$}).

Now $gh^{-1}\in \cent_G(\langle Z,Z^-\rangle) 
=L$~(cf.~(\ref{derived group})). 
Find a straight-line program of~length
$O(\log q)$ from  $\SSS^*_L$ to   $gh^{-1} $ using
Theorem~\ref{previous algorithms}(iii).  
This produces the desired   straight-line program to $g$.

\rel
$\ge 1-1/2^8$: we obtain $y$ with probability $>1-1/2^{10}$ by Lemma~\ref{opposite long root element}, and both calls to 
Lemma~\ref{transitivity of $Q$} succeed with probability $>1-2/2^{10}$.   
(N.\:B.--Recall that we are assuming that $G_0=G$, in which case 
 Corollary~\ref{random G0}
provides us with a random element of $G$ and hence a random conjugate $Z^{gy}$ of $Z$.   We will test this assumption in Corollary~\ref{verification}.)

\timing
 $O(\xi q \log q+ \mu
q\log ^2\hspace{-1pt}q)$ in (iii), dominated by the  time to find the elements 
$u$  and  $ u'$ using
Lemma \ref{transitivity of $Q$}.  (N.\hspace{.5pt}B.--It also takes   
$O(\mu q)$ time to find $h$.)
\qed
   
  \Remark \rm
\label{elements of hat G}
We have assumed in (i) that our element of $\hat
G$ was given either in terms of the Bruhat decomposition or as an automorphism of the Lie algebra for $\hat G$.
In the latter situation, 
the input to the algorithm in  \cite[Theorem~8.1]{CMT} 
or \cite{CHM}  is a linear transformation and the algorithm carries out a form of row reduction to get the Bruhat form. This is essential for our use in the Appendix, 
and nicely parallels the classical group situation \cite{KS1}.
In fact, \cite{CMT,CHM} deal  with the same question for a variety of irreducible representations of $\hat G$.

Alternatively, $\hat g$ could 
just be given as a {\em word} in $\hat \SSS$.
This possibility has already been considered:  in 
\cite[pp.~44-45]{Riebeek} and \cite{CMT}  there are deterministic algorithms which, when given $g$ as a word  
in $\hat
\SSS$, uses the relations 
(\ref{relation1})--(\ref{relation4}) or
(\ref{relation4'}), together with  
(\ref{TG action})--(\ref{W action on root groups}),
in order to rewrite  $g$  as an element $unu' $  as above. 

  In (iii) an element of the black box group $G$ is given as a string, it is not necessarily given   in terms of any available generating set.  This is essential for uses of Theorem~\ref{Main Theorem}   such as  Corollary~\ref{Upgrade corollary}. 
 
\Remark  Alternative approach to \rm (i)  \emph{\hspace{-1.5pt}avoiding}  \rm  \cite{CMT,CHM,Riebeek}:
\label {Alternative approach to (i)}
Apply the algorithm in Proposition~\ref{SL P}(iii) to the given element $\hat g\in \hat G$  (this uses  Lemmas~\ref{transitivity of $Q$} and  \ref{Expressing}
for $\hat G$).  

Here $\hat g$ might once again merely be known as an automorphism of the associated Lie algebra.  This routine has the disadvantage of requiring more time and being probabilistic;  its advantage is that it uses the present paper's relatively standard black-box methodology
employed  in (iii).

\begin{corollary}
\label {verification} 
There is an $O \big(|\SSS|  \log | \SSS |(\xi q \log q+\mu q \log ^ 2\hspace{-1pt}q ) \big)$--time  
Las Vegas algorithm which$,$ with probability $\ge 1-1/2^6,$
 checks that  $G=G_0$.
\end{corollary}

\proof
Recall  that   $G$ is given as $\langle \SSS\rangle$.  In order to prove  that $\Psi $ is an   epimorphism  we   verify  that every generator    $s\in \SSS$ lies in
$G_0$ by applying  Proposition~\ref{SL P}(iii)  to each   $s$
up to $\lceil  \log | \SSS | \rceil $ times.

\rel $ \ge 1-1/2^6$: the applications of 
 Proposition~\ref{SL P}(iii) for a single   $s\in \SSS$ all fail   with
probability $< 1/2^{7\log |\SSS|} \le 1/(2^6|\SSS|)$,
so that at least one of our  tests fails for some $s\in \SSS $  
with probability $< |\SSS|\hspace{-1pt}
   \cdot\hspace{-1pt}1/(2^{6 }|\SSS|) $.

\timing $O \big(|\SSS|\hspace{-.6pt}  \log | \SSS |(\xi q \log q+\mu q \log ^ 2\hspace{-1pt}q ) \big)\hspace{-.6pt}$ using  
Proposition~\ref{SL P}(iii) to obtain  
straight-line programs  from  $\SSS^*$  to each   $s\in \SSS$.
 \qed

\medskip
The timing in the preceding result differs from \cite[p.~145]{KS1} since the  membership test used there is deterministic, unlike our 
Proposition~\ref{SL P}(iii).

%  \hspace{.6pt}

%\newpage
\subsection {Proof of Theorem~\ref{Main Theorem} for rank $>2$}
\label{unknown type}
\label {Theorem 1.1 for rank $>2$}

In Section~\ref{Finding and verifying the presentation} we  produced a
homomorphism $\Psi \colon \! \hat G\to G$  with image $G_0$.
We  consider the various parts of Theorem~\ref{main
theorem}.%
\smallskip 

(i)  We already used \cite{BKPS}.
\smallskip 

(ii)  See 
Sections~\ref {Finding and verifying the presentation}
and
\ref{Straight-line programs}.

\smallskip 
(iii,iv,vii)
 See  Proposition~\ref{SL P}. 
\smallskip 

(v) This follows from Theorem~\ref{previous algorithms}(i)
in view of 
the new generators $X_\a(f_k)$ we  introduced
in Sections~\ref {Finding and verifying the presentation}
and
\ref{Straight-line programs}.

\smallskip 
 (vi) The second part is Corollary~\ref{verification}.
 
 The first part is the content of Sections~\ref {Start of
proof}--\ref{Finding and verifying the presentation}.
The probability of success  is at least $1/2$,
and the total time is as stated, 
 due to  all of the individual probabilities and times obtained earlier.   
 \smallskip
 
(viii)  Find $\zent(\hat G)$ using Proposition~\ref{SL P}(ii), and then 
find $\zent( G)=\zent(\hat G)\Psi $ using
Proposition~\ref{SL P}(i).  
 \qed

%\newpage

\section {Rank $2$ groups}
\label{Rank $2$ groups}

We now turn to the groups $ G_2(q) $  and $\D_4(q)$.
For the most
part we will be able to mimic and  simplify  the previous approach.
However, there are differences, such as the  use of a subgroup 
$L$ that  does not contain any long root elements.

We assume that $q>9$ in order to  avoid some exceptional situations.
In particular, we will always have $\hat G\cong G$ \cite[p.~313]{GLS},
where $\hat G$ will be known and ``concrete'' whereas $G$ will be a black box group.

\subsection{Background}
\label{Background rank 2}
In addition to $\F=\F_{q}$  we need to consider
$\F'=\F_{q^\epsilon}$, where $\epsilon$ is 1 for $G_2(q) $  and 3
for   ${^3\!D_4}(q)$. 
We retain our notation from Section 2, except that now  
$\F'$ is  $\F_q$   or $\F_{q^\epsilon}$ and 
$\{f_1,\dots,f_{\epsilon e} \}$ is
an $\F_p$-basis of $\F_{q^\epsilon}$.   

\para{Presentation.}
The groups $ G_2(q) $  and $\D_4(q)$ have a root system $\Phi$ of
type~$G_2$.

First consider  $\hat G={\D_4(q)}$. 
We start with generators $x_\a(t)$,
\vspace{1pt}
where either $\a$ is long and $t\in \F$, or $\a$ is short and $t\in
 \F' = \F_{q^3}$. Define   $\Tr\colon \!\F' \to  \F$  by  $\Tr(t)
= t + t^q + t^{q^2}\!.$ Then the Steinberg relations  \cite{St} become
(\ref{relation1})--(\ref{relation3}), where the field elements are in
$\F$  or $\F'$ for $\a$ long or short, respectively, together with 
\begin{equation}
\label{3D4}\hspace{.6pt}
\begin{array}{l}
\hspace{-20pt}
  \makebox[1.38in][l]{$ [ \hat X_\alpha   (f_k), \hat X_\beta  (f_l) ] =   $}
\qquad \hspace{-10pt} \mbox{\em for} \qquad\qquad
 \vspace{4pt}
 \\ 
    \makebox[1.38in][l]{$ 1$}
\mbox{$\alpha  +\beta   \notin \Phi$}
 \\ 
  \makebox[1.38in][l]{$\hat X_{ \a+ \b} ( \epsilon_{ \a \b}  \!  f_k
f_l)$}
\mbox{\em $ \a, \b, \a+  \b$ long}  
\\
    \makebox[1.38in][l]{$ \hat X_{\alpha +\beta }\big ( \epsilon_{\alpha
\beta } \Tr(f_k f_l)\big  )$}
\mbox{\em $\alpha , \beta $ short,
$\alpha +\beta $ long} 
\vspace{2pt} 
\\
 \hat X_{\alpha +\beta } \big ( \epsilon_{\alpha \beta } (f_k^q
f_l^{q^2} + f_k^{q^2}  \! f_l^q) \big ) \hat X_{2\alpha +\beta } \big  (
\eta_{\alpha \beta }
\Tr(f_k f_k^q f_l^{q^2}) \big ) \cdot
\vspace{4pt}
 \\
\hfill
 \hat X_{\alpha +2\beta }\big ( \delta_{\alpha \beta } \Tr(f_k f_l^q
f_l^{q^2}) \big )
\vspace{4pt}
\\ 
   \makebox[1.38in][l]{~~} \mbox{\em $\alpha , \beta , \alpha +\beta $ short, 
$2\alpha +\beta ,
\alpha +2\beta $ long} 
\vspace{2pt}
\\
 \hat X_{\alpha +\beta } ( \epsilon_{\alpha \beta } f_k f_l ) \hat 
X_{2\alpha +\beta } (
\epsilon'_{\alpha \beta } \, f_k^q f_k^{q^2} \! f_l ) \hat X_{3\alpha
+\beta }     (\epsilon''_{\hspace{-0.7pt}\alpha
\beta }  f_k f_k^q f_k^{q^2}  \! f_l )  \cdot
\vspace{4pt}
\\
\hfill  
\hat X_{3\alpha +2\beta }
(2\epsilon''{}'_{\!\!\hspace{-0.7pt}	\!\alpha \beta }
 f_k f_k^q f_k^{q^2}  \! f_l^2 )
\vspace{4pt}
\\
   \makebox[1.38in][l]{~~} 
\mbox{\em $\alpha , \alpha +\beta , 2\alpha +\beta $
short, $\beta ,  3\alpha +\beta , 3\alpha +2\beta $ long} 
\end{array}
\end{equation}
for all basis elements $f_k,f_l$ of $\F$ or $\F'$
(as appropriate).
Once again the coefficients
$\epsilon_{\a\b},\eta_{\a\b},\delta_{\a\b},$ $
\epsilon'_{\a\b},
\epsilon''_{\hspace{-0.7pt}\a\b},$ $
\epsilon''{}'_{\!\!\hspace{-0.7pt}	\!\a\b}$ are $\pm 1$ and depend only on
$\alpha$ and
$\beta$. 
Once again the right hand sides are viewed as  
products of  powers of generators $\hat X_\gamma(f_m)$ for the roots  
$\gamma$ appearing on the right side.

We again use
(\ref{Definition of $h alpha$}),
where $t\in \F'{} ^*$ when   $\a$ is short.  
Then the  analogues of 
(\ref{TG action})  and 
(\ref{W action on root groups}) 
hold.
For example:
\begin{equation}
\label{TG action rank 2}
%\hspace{-35.7pt}
\begin{array}{lll}
\hat h_{\a}(t)\hat X_{\beta}(u)\hat h_{\a}(t)^{-1} = \hat
X_{\beta}(t^{A_{\a,\beta}}u)  
&
%\hspace{1pt}
 \mbox{\em except for the next
instance}%
\vspace{3pt}
\\\hat h_{\a}(t)\hat X_{\beta}(u)\hat h_{\a}(t)^{-1} = \hat
X_{\beta}((t  t ^qt ^{q^2}) ^{A_{\a,\b}/3}u)   
&
%\hspace{1pt} 
\mbox {\em $\alpha\,$  
short,
$\beta$   long.}%
\end{array}
%\hspace{7pt}
\end{equation}

For $G_2(q)$ we obtain the required presentation 
by restricting all of the above field elements to $\F$.

We include a sketch of a proof of the second line in (\ref{TG action rank 2})  when $\hat G={^3\!D_4}(q)$.
 The twisted root system for 
  ${^3\!D_4}(q)$  has a base $\{\a,\b\}$ 
 arising from a base $\{\a_1, \a_2, \a_3, \a_4\}$ of a $D_4$-root system, where $\b=\a_2$ is the central node  and 
 $\a$ corresponds to $\{\a_1,\a_3,\a_4\}.$
We will follow \cite[pp.~233-237]{Ca1}.
If  $u\in \F$  and $t\in \F_{q^3}$  then  
$\hat X_\b(u)= \hat X_{\a_2}(u)$ and
$\hat h_\a(t)= \hat h_{\a_1}(t)  \hat h_{\a_3}(t^q)
   \hat h_{\a_4}(t^{q^2})$.  Moreover,
 
 \smallskip
 \hspace{-12pt}
$
\begin{array}{lllll}
\hat h_\a(t)\hat X_\b(u)\hat h_\a(t)^{-1}  &
\hspace{-8pt}= \hspace{-1pt}
\hat h_{\a_1}(t)  \hat h_{\a_2}(t^q)  \hat h_{\a_3}(t^{q^2})
 \hat X_{\a_2}(u)
\hat h_{\a_1}(t)^{-1}  
 \hat h_{\a_2}(t^q)^{-1}  
\hat h_{\a_3}(t^{q^2})^{-1}
 \smallskip
 \\
&\hspace{-8pt}=\hspace{-1pt}
\hat X_\b\big (t^{A_{\a_1, \a_2}} (t^{q})^{A_{\a_3, \a_2}} (t^{q^2})^{A_{\a_4,\a_2}} u\big )
 \smallskip
\end{array}$

\noindent
with  $A_{\a_1,\epsilon \a_2}=
A_{\a_3,\epsilon \a_2}= A_{\a_4,\epsilon \a_2}
= \epsilon = A_{\a, \epsilon\b}/3 $ 
for $\epsilon=\pm1$,
which implies the second assertion in (\ref{TG action rank 2}).

\para{The subgroup $ \hat S$.}
For both $G_2(q)$ and ${^3\!D_4}(q)$ the subgroup $\hat S$ generated by
 the long root groups $\hat X_\a$ is isomorphic to $\SL(3,q)$.

\para{The subgroups $ \hat Q$ and $ \hat L$.}
If $ \hat Z$ is a long root 
subgroup of $ \hat G$ and  $1 \neq z \in  \hat Z$, then 
$\cent_{ \hat G}(z) =
\cent_{ \hat G}(\hat Z) = \hat  Q \semi  \hat L$ with
$ \hat Q$ and $ \hat L$ as follows:
$$
\begin{array}{|l|l|l|l|l|c|l}
    \hline 
 \hat G~~~
 & G_2(q), q \neq 3^a ~~~& G_2(q), q= 3^a~~~ & \D_4(q) 
\raisebox{-1ex} {~}\raisebox{2.5ex} {~}
\\
    \hline 
 \hat L
 & \SL  (2,q)& \SL  (2,q)& \SL  (2,q^3)
 \raisebox{2.5ex} {~}
\\
    \hline \hat Q
  &q^{1+4} & q^{1+(2+2)} & q^{1+8} 
   \raisebox{2.5ex} {~}
\\
    \hline
 \hat T
  &\Z_{q-1}\times \Z_{q-1} & \Z_{q-1}\times \Z_{q-1} & \Z_{q-1}\times
\Z_{q^3-1}
 \raisebox{-1.3ex} {~}\raisebox{2.5ex} {~}
\\
    \hline
T_{ \hat L}
  &\Z_{q-1}  & \Z_{q-1}  & \Z_{q^3-1}
  \raisebox{-1.2ex} {~} \raisebox{2.2ex} {~}
\\
    \hline 
\end{array}
$$
where we have   included the structure of maximal tori $\hat T$ of $\hat G$
and
$T_{\hat L}$ of $\hat  L$.

\begin{lemma}
\label {rank 2 ppd}
\begin{itemize}
\item[\rm(i)]
With probability $\ge 1/3 q,$ an element $\tau\in
G_2(q)$ has order  $p\hspace{-1pt}
   \cdot\hspace{-1pt} \ppd^\#(p;2e);$ and then $\tau^{q+1}$
is a long or short root element.
\item[\rm(ii)]
With probability $\ge 1/9 q,$ an element $\tau\in
{^3\!D_4}(q)$ has order  $p\hspace{-1pt}
   \cdot\hspace{-1pt} \ppd^\#(p;6e),$ and
then
$\tau^{q^3+1}$ is a long root element.  

\end{itemize}
\end{lemma}

\proof  We first construct elements of the indicated orders. 
There is  a central product  $\SL  (2,\F')\circ \SL  (2,q)$ 
of  a short root $\SL  (2,\F')$ and
  a long root $\SL  (2,q)$, and this contains  elements of the desired
order. 
As  in Lemma~\ref{tori
probability}, an element $\tau$ of the stated order lies in a parabolic,
hence in a central product as above, and hence powers to a root element.

The probability estimates are  
obtained as in Lemma~\ref{tori
probability}, but   
are   simpler.~\qed

\medskip
{\em Opposite} long root elements and root groups are defined as in
Section~\ref{Properties of $G$}.

% \newpage

\begin{lemma}
\label {SL 3 test}
Let $z$ be a long root element.
\begin{itemize}
\item[\rm(i)]
{\rm(\ref{opposite test})}  holds.
\item[\rm(ii)]
{\rm Lemma~\ref{opposite long root element}} holds.
\item[\rm(iii)]
{\rm Lemma~\ref{S J probabilities2}(i)} holds.

\item[\rm(iv)]
All long subgroups isomorphic to  $\SL  (3,q)$
are conjugate.
\item[\rm(v)]
If $p\ne 3$ then three short root elements of $G_2(q)$ never generate a
group isomorphic to  $\SL  (3,q)$.
%  or xxxxxxxxxxxxxxx$\PSL(3,q)$. 

\end{itemize}
\end{lemma}

\proof
(i) This  follows from the analogue of Lemma~\ref{commutator
relations}.

(ii,iii) These are proved exactly as in Section~\ref{Probability and long
root elements} (cf. Table~\ref{n(S,R)}). 

(iv) See
 \cite{Cooperstein} or   \cite{Kl1,Kl2}.

(v)
See  \cite{Kl1}.  \qed
\medskip 

Of course, the conclusion in   (v)   is false for $p=3$
due to  the graph  automorphism of $G_2(q)$.

\subsection {Finding a root group $Z$ and the subgroups 
$Z^-$, $R$ and  $S$}
\label{Root element and $J$}
\label{Root element and S}
As in Section~\ref {Start of proof}, 
we now consider a black box group $G$ 
that is a nontrivial homomorphic image of the
universal cover $\hat G$   of $G_2(q)$ or
$\D_4(q)$.
Since  $q>4$,   $ \hat G \cong G\,$
\cite[p.~313]{GLS}.
Find the probable type of $G$    using 
\cite{BKPS}. 

%\newpage

 We now imitate parts of Sections~\ref{Finding a root element}  and \ref{Finding a root group}.  
Choose up to $ 90 q$  elements $\tau$ in order to find one of order 
$pl \! = \! p\hspace{-1pt}
   \cdot\hspace{-1pt} \ppd^\#(p;2\epsilon e)$;
for $z:=\tau^{q^{\e}+1}$ choose up to $ 120 $ pairs $z',y$   of conjugates of  $z$;
for each pair, test whether both are opposite $z$ and whether
 $S:=\langle z,z',y\rangle$  and
$ S_2:=\langle z, z' {}^{\tau ^p },y
\rangle$ are both isomorphic to ${\hat S}=\SL(3,q)$;
\vspace{.5pt}
 and, if so,  find constructive isomorphisms
$\Psi _S\colon    \!{\hat S}  \to S   $  and 
$\Psi_{S_2} \hspace{-1pt}\colon \!{\hat S} \to S_2$, together with generating sets $\SSS_{\hat S}$  and $\SSS_S^*$  of $\hat S $ and $S   $, respectively, such that $  \SSS_{\hat S}\Psi_S=\SSS_S^*$.  

Find the root groups 
$Z  $  and  $Z^{-}$ in $S$ such that $z\in Z$ and $ z' \in Z'$.

Let  $R: =\langle Z,Z^{-}\rangle \cong  \SL (2,q)$.

As in Section~\ref{Finding a root group}, find  
$c\in \Oh_p\big (\cent_{S_2}(Z)\big)  \Oh_p\big(\cent_S(Z)\big) \subseteq
\Oh_p\big(\cent_G(Z)\big)$    such that  
 {\em $\tau^pc$ centralizes $  R$ and has order
divisible by $l$.}

{\em Correctness}: 
By Lemma~\ref{rank 2 ppd}(ii),
 the element $z$ just constructed  is a long root element 
if $G$ is ${}\D_4(q)$.  
If $G$ is $G_2(q)$ and $p=3$,  
it makes no difference  whether we are using long
or short root elements, since these are conjugate in $\Aut G$,  so we
may assume that $z$ is long.  
If $G$ is $G_2(q)$ and $p\ne3$ then we might have obtained a short root element $z$, but then we will not obtain $S\cong \SL(3,q)$
 by Lemma~\ref {SL 3 test}(v).

\rel $\ge 1-1/2^{9}$.  
 For, a choice  $\tau$ has the correct order and produces a long root element with probability
  $\ge 1/9q$  by Lemma~\ref{rank 2 ppd}(i), so that 
  we fail to obtain an element $\tau$ of the desired type with probability $\le (1-1/9q)^{90q}<1/2^{10}$.
The tests involving  a single choice $z',$ $y,$ $S,$ $S_2$ all succeed with probability $\ge(1/12)(1/3)(1/4)^2(1/2)^2>1/2^{12} $ (by Lemma~\ref {SL 3 test}(ii,iii) and Theorem~\ref{previous algorithms}(ii)),
 so that  
 the tests for all $120$ pairs $z',y$
   all fail with probability $<(1-1/2^{12} )^{120 }<1/2^{10}$.

\timing $O(\xi qe + \mu q\log^2\hspace{-1pt}q)$, dominated by 
finding $\Psi _S$  and $\Psi _{ S_2}$.

\subsection{The subgroups   $L$, $T$    and $N$.}
\label{The subgroups $L$ and $Q$}

We will use additional subgroups analogous to ones in 
Sections~\ref {Finding a root element}--\ref{The root groups
$Xalpha$}.

\para{The group\,}$L:=\langle \cent_S(R), \tau^p c \rangle $
is a subgroup of $\cent_G(R)=\SL(2,q^\epsilon)$ of order divisible by
both 
$|\cent_S(R)|=q-1$ and $|\tau^p c|$, which is a $\ppd^\#(p; 2 \epsilon
e)$.  Then  $L=\cent_G(R)$  since $\SL(2,q^\epsilon)$  has  no such proper
subgroup  for $q>9$ \cite[Sec.~260]{Di}.   

As in Lemma~\ref{one class of SL 3}(ii), {\em the pair~$(S,L)\!$ is uniquely 
determined up to conjugacy~in~$G$.}%
 
Use Theorem~\ref{previous algorithms}(ii) up to 10 times in order to obtain a constructive 
isomorphism $\Psi _L\colon \! L\to\SL  (2,\F')$.

\rel $\ge 1-1/2^{10}$.

\timing $O(\xi |\F'|\log q + \mu |\F'| \log^2 \hspace{-1pt}q)$
to obtain $\Psi _L$.
%\medskip 

\para{The subgroups $T_S$, $T_L$, $T$ and $N$.}
First note that $ \hat G$ acts transitively by conjugation on the set of triples $( \hat S_1, \hat R_1,T_{\hat S_1})$  with $T_{\hat S_1}$ a maximal split torus of $ \hat S_1\in  \hat S^{ \hat G}$ normalizing
${\hat R_1}\in  \hat R^{ \hat G}$.  Hence $ \hat G$ is also transitive  on the set of 
4-tuples $( \hat L_1,\hat S_1,T_{\hat L_1}, T_{\hat S_1})$
  with $T_{\hat S_1}$ a maximal split torus of $ \hat S_1$
normalizing  
$ \hat L_1=\cent_{\hat G}(\hat R_1)$ and centralizing a 
(unique) maximal split torus $T_{\hat L_1}$
of $\hat L_1$ (which must therefore contain the torus
 $\cent_{\hat S_1}(\hat R_1)= \hat S_1\cap  \hat L_1$). 
Then
$T_{\hat S_1}T_{\hat L_1}$ is a maximal torus of
$ \hat G$  and is normal in 
$\langle  \norm_{\hat S_1}(T_{\hat S_1}), \norm_{\hat L_1}(T_{\hat L_1}) \rangle$ 
\vspace{2pt}
 (compare Lemma~\ref {L, S and tori}).

With this in mind, use $\Psi_S$ (and Theorem~\ref{previous
algorithms}(iii,iv)) to find a maximal split torus $T_S$ of $S$
normalizing $R$,  $Z$ and $Z^-$. Then
$T_S$ normalizes $\cent_G(R)=L$,   and hence normalizes and so  centralizes the
unique  maximal split torus $T_L\ge S\cap L$ of $L$
(by the preceding paragraph).  Find $T_L$   
using
$\Psi_L$.  
(Compare  Lemma~\ref{matching}  --  but here we are only
working with a 2--dimensional vector space.  Moreover, unlike in the
large rank case,  the torus $T_L$ is uniquely determined
by $S\cap L$.)

 Then $T:=T_ST_L$ is a maximal 
torus of
$G$.

Find $\norm_S(T_S)$
and $\norm_L(T_L)$  using  $\Psi_S$
and $\Psi_L$.  The above observations concerning $\hat G$
imply that   $T \nor N:=\langle  \norm_S(T_S), \norm_L(T_L) \rangle $ 
and $N/T$ is the Weyl group of $G$.  

From this point on we  will no longer explicitly use $S$.

 \vspace{-2pt}

\subsection{Root groups. }
\label{Root groups.} 
 ~Let $\{\alpha_1,\alpha_2\}$ be a base for 
 a root system  $\Phi$ associated with the Weyl group
 $N/T$, with $\alpha_1$
long.  Let
$\nu= 2\a_1+3\a_2 $ be the highest root, and label $Z=X_\nu$ and
$Z^-=X_{-\nu}$.

Find the two (short!)  root groups  
of $L$ normalized by $T$ using $\Psi_L$,
pick one of them and label it $X_{\a_2}$, 
then label the other one $X_{-\a_2}$. 
The $N$-conjugates of $X_{\pm \nu}$ and
$X_{\pm\a_2}$ are the 12 root groups of $G$  normalized   by $T$; %
the action of $N$ labels each as   $X_\alpha$ with 
$\alpha\in \Phi$.

Coordinatize $L$ using $\Psi_L$, obtaining 
 $X_{\pm\a_2}(f_k)$, 
$n_{\a_2}(1)$  and 
$  h_{ \a_2}(f_k)  $  for    $f_k\in \F'$.

\timing $O(\mu  q \log ^2q )$, dominated by 
$O(e)$ uses of  Theorem~\ref{previous algorithms}(iii,iv) for
 $\Psi_ L$.
 
\medskip

As in Section~\ref{Finding and verifying the presentation}, 
we next  show that $\hat G$  maps onto  $G_0:=\langle  X_\a\mid \a\in
\Phi\rangle$:

\begin{proposition}
\label{constructing lambda rank 2}
\  ~ There is a deterministic  $O(\mu  \log^2 \hspace{-1pt}q)$-time    algorithm 
that  labels any given element of any root group
$X_\a, $  $\a\in \Phi,$  as $X_\a (t)$ for some $t\in \F$ or
$\F',$ in such a
way  that the map $ \hat X_\a (f_k)\mapsto X_\a (f_k)$ $($for all
appropriate $\alpha $  and $k)$ extends to an epimorphism
$\Psi\colon \! \hat G\to G_0 $.%

\end{proposition}

\proof
By  (\ref{TG action rank 2}), $\hat T_L$ acts 
transitively  on
the nontrivial elements of $\hat X_{\a_1}\hspace{-1pt}.$  Thus, we can choose any nontrivial element  
$ X_{\a_1}$ and label  it
$X_{\a_1}(1)$,~after which  the remaining   
labels $X_{\a_1}(f_k)$
are forced    by  (\ref{TG action rank 2}). Namely,    
$h_{ \a_2}(t)^{-{A_{\a_2,\a_1}/3}}$
conjugates $X_{\a_1}(1)$ to $X_{\a_1}(tt^{q}t^{q^2})$.
Applying this for distinct $t=f_k, f_k+1, af_k+1$ in $\F$ gives us  $X_{\a_1}(u)$ for  $u=f_k^3, (f_k+1)^3 $ and $ (af_k+1)^3$.   We may assume that $p\ne3 $ (as otherwise the elements $f_k^3$ span $\F$).
Since $f_1=1$, it is easy to see that we now have obtained
all of the elements $X_{\a_1}(f_k)$.

 By the rank 2 analogues of (\ref{Definition of 
$h alpha$}) and (\ref{W action on root groups}), we can now coordinatize  
$X_{\a_1}^{n_{\a_2}(1)}=X_{-\a_1-3\a_2}$. 

By   (\ref{3D4}), $[ [X_{ \a_2 }(1),X_{ \a_ 1}(f_k)], X_{ \a_1}(1)]\!=  \!     X_{ 2\a_1+3\a_2}(\epsilon_{ \a_1+3\a_2,  \a_1 }^{~} \epsilon_{
 \a_2, \a_1}''f_k)$ whenever $ f_k\in \F$, so we can coordinatize
  $X_{ 2\a_1+3\a_2}$.

For each element $x$ in an $\F_p$-basis of   $X_{ -\a_ 1}$,
find its coordinate $u$ via the relation  
$[ X_{ 2\a_1+3\a_2 }(1),x ] =
X_{ \a_1+3\a_2 }(\epsilon_{ 2\a_1+3\a_2, - \a_1 }u) $
in (\ref{3D4})  
by  using linear algebra in $X_{ \a_1+3\a_2 }$.
This produces the coordinates of a basis of $X_{ -\a_ 1 }$
and hence of any given element of $X_ { -\a_ 1 }$.

Use (\ref{Definition of 
$h alpha$}) and (\ref{W action on root groups}) to   coordinatize all $\langle n_{\a_1}(1) ,  n_{\a_2}(1) \rangle$-conjugates of $X_{\a_1}$ and  $X_{\a_2}$, and hence of all root groups $X_\a$.

Thus, we have obtained a   map  $\Psi\colon  \!\hat X_\alpha(f_k)\mapsto
X_\alpha(f_k)$  (for all appropriate  $\alpha$ and~$k$).
Verify (\ref{3D4}) in order to show that  $\Psi$  extends to
an epimorphism 
$\hat G\to G _0$. 
As in the proof of Proposition~\ref{constructing lambda}, this algorithm
is deterministic, and runs in the stated time.
\qed

 \vspace{-2pt}

\subsection {Linear algebra in $Q/Z$}
\label {Linear algebra in $Q/Z$ rank 2 case}

Next we imitate Section~\ref{Linear algebra in $Q/Z$}.

 \vspace{-4pt}
\para{Effective transitivity  of the subgroup $Q$.}     
  Lemma~\ref{transitivity of $Q$} {\em holds for}  $Q:=\langle 
X_\alpha\mid  \alpha\in  \Phi^+
\rangle$, 
 using the  exact same proof, still requiring
$O(\xi q\log q + \mu q\log^2 \hspace{-1pt}q)$ time and still succeeding 
with probability $>1-1/2^{10}$.

 \vspace{-4pt}

\para {Linear algebra in $Q/Z$.}
If we exclude $G_2(q)$ with $p=3$, this is the same as in
Proposition~\ref{expressing}.  Namely, $Q$ is still of ``extraspecial
type"  (i.\hspace{1pt}e., it behaves exactly as in Lemma~\ref{pairing roots}(ii)), and
we can again peel off the root elements by commutations
as in the proof of Proposition~\ref{expressing}.

However, since this ``peeling'' involves traces of field elements, we will
be more careful.
 List the positive roots $2{\a_1}+3{\a_2} =
\nu,$ ${\a_1}+3{\a_2},$ ${\a_1}+2{\a_2},$ ${\a_1}+{\a_2},$
${\a_2},$ ${\a_1}$. Our given  $g \in Q$ can be written  $g =
\prod _{\gamma\in
\Phi^+
\setminus
\Phi_L} X_\gamma(t_\gamma)$ in this order, and we must find the field elements
$t_\gamma$.

By   (\ref{3D4}), ${ X}_{{\a_1}+3{\a_2}}$ commutes with  the positive root groups other
than ${ X}_{\a_1}$. 
Since
$[g,{ X}_{{\a_1}+3{\a_2}}(1)] = 
X_{2{\a_1}+3{\a_2}}(\epsilon_{{\a_1},{\a_1}+3{\a_2}^{~}}\!t_{ {\a_1} })
$, 
as in Lemma~\ref{pairing roots}(i)
we   deduce $t_{{\a_1}}$  using linear algebra in $\F$.

Let $g_1:=g X_{{\a_1}}(t_{{\a_1}})^{-1}$.
By   (\ref{3D4}),
$$
\begin{array}{llllllll}
g_1'& \hspace{-7pt} :=[g_1,{ X}_{{\a_1}+{\a_2}}(1)]=
[X_{\a_1+2\a_2}(t_{\a_1+2\a_2}),X_{\a_1 + \a_2}(1)][{ X}_{{\a_2}}(t_{{\a_2}}),{
X}_{  {\a_1}+{\a_2}}(1) ] 
 \vspace{2pt}
\\
 & \hspace{-4.2pt} =  { X}_{{\a_1}+2{\a_2}} \big (\epsilon_{{\a_2},
{\a_1}+{\a_2}}^{~}(t_{{\a_2}}^{^{\scriptstyle q}}+t_{{\a_2}}^{^{\scriptstyle q^2}}) \big )  {
X}_{{\a_1}+3{\a_2}} \big (\eta_{{\a_2}, {\a_1}+  {\a_2}}^{~}
\Tr(t_{{\a_2}}   t_{{\a_2}}^{^{\scriptstyle q}}   )  \big )\cdot
 \vspace{2pt}
\\
&\hfill  {X}_{2{\a_1}+3{\a_2}  } \big (\delta_{{\a_2},
{\a_1}+{\a_2}}\Tr(t_{{\a_2}} ) \big )  X_{2\a_1+3\a_2} \big ( \epsilon_{\a_1+\a_2, \a_1+2\a_2}T(t_{\a_1+2\a_2})\big ) .
% \vspace{-2pt}
\end{array}
$$
Then $[g_1',X_{\a_1}(\epsilon_{{\a_1+3\a_2},{\a_1}} 1 )] 
= X_{{2\a_1}+3{\a_2}} \big(\eta_{{\a_2}, {\a_1}+  {\a_2}}^{~}
\Tr(t_{{\a_2}}   t_{{\a_2}}^{^{\scriptstyle q}}   )  \big)$ 
gives us
$\Tr(t_{{\a_2}}   t_{{\a_2}}^{^{\scriptstyle q}}   )$.

Also,%
\smallskip

 \hspace{-19pt}
$
\begin{array}{llllllll}
%\hspace{-4.2pt}
[g_1',X_{-\a_1}(\epsilon_{{\a_1+3\a_2},{-\a_1}} 1 )] 
= X_{{\a_1}+3{\a_2}}\big(\delta_{{\a_2},
{\a_1}+{\a_2}}\Tr(t_{{\a_2}}) + \epsilon_{\a_1+\a_2, \a_1+2\a_2}T(t_{\a_1+2\a_2})\big),
%\end{array}
\smallskip
\\
%\begin{array}{llllllll}
g_1'[g_1',X_{\a_1}(\epsilon_{{\a_1+3\a_2},{\a_1}} 1 )]^{-1}[g_1',X_{-\a_1}(\epsilon_{{\a_1+3\a_2},{-\a_1}} 1 )]^{-1} 
\\
%&
 \hspace{110pt}
=
{ X}_{{\a_1}+2{\a_2}} \big (\epsilon_{{\a_2},
{\a_1}+{\a_2}}^{~}(t_{{\a_2}}^{^{\scriptstyle q}}+t_{{\a_2}}^{^{\scriptstyle q^2}}) \big ).
\end{array}
$

\noindent
Hence,   we deduce $t_{{\a_2}}^{^{\scriptstyle q}}+t_{{\a_2}}^{^{\scriptstyle q^2}}$.
The identity $ (t_{{\a_2}}^{^{\scriptstyle q}}+t_{{\a_2}}^{^{\scriptstyle q^2}})
(t_{{\a_2}}^{^{\scriptstyle q}}+t_{{\a_2}}^{^{\scriptstyle q^2}})^{^{\scriptstyle q}} = (t_{{\a_2}}^{^{\scriptstyle q}})^2 
+ \Tr(t_{{\a_2}}   t_{{\a_2}}^{^{\scriptstyle q}}  ) $
along with  $\Tr(t_{{\a_2}}   t_{{\a_2}}^{^{\scriptstyle q}}  ) $
give  us  $(t_{{\a_2}}^{^{\scriptstyle q}})^2$ and hence also
$  t_{\a_2} ^{^{\scriptstyle q}} $ up to sign. Since we already know
 $  t_{{\a_2}}^{^{\scriptstyle q}}+t_{{\a_2}}^{^{\scriptstyle q^2}} $, we  
deduce $  t_{\a_2} ^{^{\scriptstyle q}} $ and hence also  $t_{{\a_2}}$.

The same procedure, with  the roles of $\a_2$ and $\a_1 + \a_2$ reversed, 
yields $t_{{\a_1 + \a_2}}$.

Let $g_2:=g_1X_{{\a_2}}(t_{{\a_2}})^{-1}X_{{\a_1+\a_2}}(t_{{\a_1+\a_2}})^{-1}$. 
\vspace{1pt}
As above,
$$[g_2, X_{{\a_1}}(\epsilon_{{\a_1}+3{\a_2},{\a_1}}1)] = X_{2{\a_1}+3{\a_2}}(\epsilon_{{\a_1}+3{\a_2},{\a_1}}^{~}
t_{{\a_1}+3{\a_2}^{\,\!\!}} ) $$  
 yields 
$t_{{\a_1}+3{\a_2}}$.    We obtain  $ t_{2{\a_1}+3{\a_2}}$ and
 $t_{{\a_1}+2{\a_2} }$ similarly.
 
As in Proposition~\ref{expressing}, this linear algebra routine is
deterministic, and takes   
$O(\mu \log q)$ time.

\subsection{Proof of Theorem~\ref{main theorem} for rank $2$}
\label{Proof of Theorem}
We can now complete the proof of  Theorem~\ref{main theorem}.

\para{Straight-line programs.} The analogue of 
 Proposition~\ref{SL P} is proved in the
same manner as in that proposition.  The timing 
for the analogue of Proposition~\ref{SL P}(iii) is
 $O(\xi |\F'|\log q 
+ \mu |\F'|\log q)$, dominated by finding the elements $u$
and $u'$ occurring  in the proof  of Proposition~\ref{SL P} and finding   straight-line  
programs in
$L$.

\para{Completion of proof.}  
This is exactly as  in Section~\ref {unknown type},
in view of Proposition~\ref{constructing lambda rank 2} 
and the analogue of 
 Proposition~\ref{SL P}.  As usual,
  (viii) is unnecessary
since     ${\zent(\hat G)=1}$.~\qed 

%\newpage
 
\section {Concluding remarks}

\label {Concluding remarks}
\label{Additional remarks} 
%\vspace{-4pt}

\begin{comments}
\label{Small $q$} 

{Small $q$.} 
\label{small q}\rm
When $q\le9$,  in place of Lemmas~\ref{tori probability} or 
\ref {rank 2 ppd}  we can simply find exact orders of elements
(replacing the stated $l$ by  $l:=|\hat T^*|$ using (\ref{table of L})).  
We still need the fact that  $q>3$ in order to have
elements behaving as in the conclusions of 
Lemmas~\ref{tori probability}(i) or
\ref{rank 2 ppd}(i). 
We also used the fact that  $q>4$ in order to avoid exceptional universal covers.  
  
When $q=9$, two opposite long root elements never generate an 
$\SL(2,9)$, but
instead generate $\SL(2,5)$.  However, 
as in Lemma~\ref{S J probabilities2}(i), inclusion of a third long root element
generates $\SL(3,9)$ with high probability, after which theç rest of our algorithm goes through.

For $q>4$, in rank $>2$  the only other needed change is (possibly) to select more elements  in order to handle the fact that the probabilities in   situations such as  
Lemmas~\ref{tori probability} and
 \ref{opposite long root element}--\ref{S J probabilities2}
% and \ref {SL 3 test}
 are no longer as high as in those lemmas.
 
However, for rank 2 a different approach is needed when $q$ is 5 or 7:
 in Section~\ref{The subgroups $L$ and $Q$}, elements of $\cent_G(R)$ of order $q-1$ and $q+1$ need not generate  $\cent_G(R)$. 
 One way is to use the fact that elements of the stated orders generate $\cent_G(R)$ with probability $>1/2$, while another 
 proceeds as in Remark~\ref{Involution centralizers.} below.

\end{comments}

\begin{comments}{Speculations on implementation.}  \rm
We expect that versions of the algorithms will be implemented.  
For rank $>2$ we suspect that there is no need to find $J$.  Instead, 
$\langle\cent_S(L) ,  \tau  \rangle$  
or $\langle\cent_S(L) ,  \tau , 
\tau _0\rangle $ appears to be the desired group $L$ when $q>2$
(in the notation of Section~\ref{Finding  $L$}).
For example, if $G$ does not have type $E_8$ then $L$ is essentially a classical group, and  
$\langle\cent_S(L) ,  \tau  \rangle$  
or $\langle\cent_S(L) ,  \tau , 
\tau _0  \rangle $ 
 acts irreducibly and primitively on its natural module.  
  Now the ppd-orders and  \cite{GPPS} can be used to obtain a small list of possibilities to check, and presumably to rule out  most of them 
by  careful examination of the  elements $\tau$ and $\tau_0$.
 
 \end{comments}

\begin{comments} {The omitted groups $\FF_4(q)$.} \rm
We expect that the groups $\FF_4(q)$ will eventually be handled in a manner resembling
Section~3.  However, those groups  involve
more intricate commutator relations than other groups of Lie type.

   The natural representation  of $\:\FF_4(q)$ is dealt with in \cite{Baa4}, assuming 
   the correctness of a  complicated conjecture concerning $\F_{q}=\F_{2^{2e+1}}$
   and of
   a conjecture  concerning
 the actions of elements of $\:\FF_4(q)$ on the
natural module.  Apparently this approach does not  work for other
absolutely irreducible  representations of  $\FF_4(q)$ in characteristic 2.

\vspace{8pt}

\em
  {\rm Remarks \ref{oracles}-\ref{Involution centralizers.}} concern variants of 
{\rm Theorem~\ref{main theorem}} that $($almost$)$ run in polynomial time.  However$,$ these have yet to be carefully checked
before there can be a sequel to this paper.  
 \end{comments}

\begin{comments} {The factor $q$ and oracles$:$ rank $>2$.} 
\label {The factor $q$.}
\label {oracles}
%\hspace{-4pt}
\rm  Our algorithm searched  for a long root element $ z \in G$, 
and then  $\langle z, z^g\rangle$ ($ g\in G$) is guaranteed to be  a proper subgroup of $G$. In fact,
 with high probability
$\langle z, z^g, z^h \rangle$ ($\,g,h\in G$) is a 
long root $\SL(3, q)$. 
Unfortunately, the probability of finding by random search
an element for which some power is a long root element 
is unreasonably  small for  groups defined over large fields. An
alternative strategy  is to search for 
semisimple elements  closely related to long root elements.

This was accomplished  in a number of the papers cited following 
Theorem~\ref{previous algorithms}.
More significantly, the factor $q$ in the
timing of analogues of Theorem~\ref{main theorem} was removed
  by  assuming the availability of an
$\SL(2,q)$-oracle to constructively recognize $\SL(2,q)$
as well as a Discrete Log oracle for $\F_{q}^*$, and possibly
also  for 
$\Z_{q+1}$  (cf. Section~\ref{Introduction}).
Then  suitable $p'$-elements were used to construct subgroups such as
$\SL(3,q)$, $\SU(3,q)$ or  $\Sp(4,q)$. 

 Here   we  comment on the requirements in order for
this approach to be used with   exceptional groups {\em of rank $>2$
when $q>4$}.
Find and use  an element $\tau$  of order $\ppd^\#(p; e) l$  \,or\,
$\ppd^\#(p;2e) l$   
in the notation of 
Lemma~\ref{tori probability}; such an element is 
 obtained as the product of elements of $R=\SL(2,q)$ and
$L= \cent_G(R)$.  The element $\tau$ needs to have two further
properties:  (a)   $\tau^l$ lies in a long $\SL(2,q)$,  and (b) two
conjugates of $\tau^l$ probably generate a subgroup containing long root
groups (in which case a long root group is obtained via constructive
recognition of the subgroup).  

Condition (a): 
As in Section~\ref{Finding a root element},  we obtain an element $\tau^l$ of a long $\SL(2,q)$, except perhaps in type $F_4$ where this might belong to a short $\SL(2,q)$.  In the latter case, we  obtain an element of order $p$ of this $\SL(2,q)$, and then proceed exactly as in Section~\ref{Finding a root element} to  distinguish
long and short root elements in odd characteristic 
(or use Remark~\ref{Involution centralizers.} below).

 There is a  problem with the first element order  $|\tau|=pl$ 
  in Lemma~\ref{cyclic tori} for $E_7(q)$. This is the only instance with a factor 
$  \ppd^\#(p;e )\ppd^\#(p;2e)$.
 One   way  around this difficulty is to modify Lemma~\ref{rationale for choice of tau}: 
use elements
   $\tau_1$ of order 
 $ \ppd^\#(p;e ) \ppd^\#(p;9e) $
 and $\tau_2$ of order 
 $ \ppd^\#(p;2e ) \ppd^\#(p;18e) $
 normalizing subgroups of type $E_6(q)$ and $\E_6(q)$,  respectively.   
 Once   conjugates of 
 the    powers $\tau_1^{(q^9-1 )/( q^3 - 1)}$ 
  and  $\tau_2^{(q^9+1 )/( q^3 + 1)}$ have been
arranged 
(by conjugation) to  generate a subgroup $\Spin_4^+(q)$ 
we will have two commuting long subgroups $R\cong R_1\cong \SL(2,q)$; 
and once other conjugates  of
 $\tau_1^{(q^9-1 )/( q^3 - 1)}$ 
  and  $\tau_2^{(q^9+1 )/( q^3 + 1)}$
  have been
arranged  to lie in that long subgroup $R=\SL(2,q)$
then $\<R_1,\tau_1^{q-1} ,\tau_2^{q+1} \>$ will be $L=\cent_G(R)$.
 
Condition (b):  If 
$\tau^l$ lies in a long $\SL(2,q)$ then two of its conjugates lie in the
group generated by two such subgroups $\SL(2,q)$, and hence 
for rank $>2$ everything
reverts to an orthogonal group setting \cite[Proposition~3.2]{Ka2},
where the required (probable) generation was  proved in  \cite{BrK1,
BrK2}.

Starting from a long root element obtained by generating a suitable subgroup in this manner, 
and assuming the availability of suitable oracles, 
the remainder  of our  algorithm goes through.
These oracles are  
the aforementioned ones for $\SL(2,q)$ and   $\F_{q}^*$;
and, in the $\E_6 $ case, one for Discrete Logs in $\Z_{q+1}$
(cf. \cite{Br2}).

  \end{comments}

\begin{comments}{Rank $2,$ even $q$.}
\label{Even $q$.}
 \rm
 The method in the preceding Remark  also works for type $G_2$
 in characteristic 2, 
 using an element of order $ 3\ppd^\#(p; \delta e )$ when
${\hat G} = G_2(q)$ (where $\delta $  is $2$ if  $3\mid q-1$
and  $1$ if  $3\mid q+1$).   

Unfortunately, when $q$ is even $\D_4(q)$   does
not possess any class  $x^G$ of semisimple elements for which
$\langle x, x^g \rangle$ \,($g\in G$)  is a proper subgroup
with high probability.  Therefore,  our approach   in  
Section~\ref{Rank $2$ groups}  appears to be  the only option 
for these groups.
      \end{comments}

\begin{comments}{Odd $q$ and involution centralizers.}
\label{Involution centralizers.}
 \rm
 There is a different way to handle part of 
 Theorem~\ref{main theorem}  that can produce a  long
 $\SL(2,q)$  in polynomial time   when $q$ is
odd, assuming the
 availability of suitable oracles 
 as in Remark~\ref{oracles}.
With high probability, a random element has  even order and a power is
an  involution
$t$ conjugate to the one in $R$.  (There may be other involutions
encountered, but the desired conjugacy class  will  occur with high
probability.)
 Then $\cent_G(t) =R\circ L$ can be found in polynomial time with high
probability
 \cite{Bor,Bray,HLORW,PW}, after which it is easy to
find both $R$ and
$L$.    As in Remark~\ref{oracles}, given suitable oracles
 the rest of
our algorithm appears to go through.  Note that,
using this approach,   we have already obtained the crucial subgroup
$L$, and hence there is no need for the subgroup~$J$.

In rank 2, the~$^3D_4(q)$ case
appears to need  oracles to  constructively recognize $\SL(2,q^3)$
and  for Discrete Logs in $\F_{q^3}^*$.

  \end{comments}

\begin{comments} {Rank $1$ groups.} \rm 
\label {Rank 1 groups.}
An early version of this paper contained Las Vegas algorithms for handling 
rank 1 exceptional groups --  Suzuki groups 
$\Sz(q)=\!\!\,^2\hspace{-1pt}B_2(q)$ and Ree~groups~$^2G_2(q)$~-- except for timing that involved a factor $q^2$ or $q^3$, respectively, as well as use of 
 a length $O( q^3\log^2\hspace{-1pt}q)$
presentation  for~$^2G_2(q)$. 
(This result implies that, {\em in} Corollary~\ref {Upgrade corollary},
{\em there is no need to exclude $^2G_2(q)$ composition factors.})

However,  that older  approach now seems far less interesting. 
A lovely black box Las Vegas algorithm for $\Sz(q)$ is in \cite{BrB}, with  timing involving a factor $q$.  An alternative  approach  \cite{Baa1,Baa2}  deals with $\Sz(q)$ as a matrix group 
and avoids any such factor but assumes the correctness of a  complicated conjecture concerning $\F_{q}=\F_{2^{2e+1}}$.

The   Ree groups $^2G_2(q)$ were studied in \cite{Baa3}  
as 7-dimensional matrix groups
using an involution centralizer and an $\SL(2,q)$ oracle (cf. Remarks~\ref{oracles} and \ref{Involution centralizers.}), this time assuming a   complicated conjecture concerning  the field
$\F_{q}=\F_{3^{2e+1}}$. 
  There is some hope that a different use of an involution centralizer (together with suitable oracles) can handle the black box setting without a need for any such conjecture or any factor $q$ in the timing. 

 \end{comments}
 %\newpage
 
 \noindent
\emph{Acknowledgement}:  We are grateful to
the referee, as well as to  Peter Brooksbank and \'Akos Seress,
for many helpful comments
and suggestions. 
 
% indebted to
\appendix

\section*{Appendix: The group $\hat E_7(q)$
and its Lie algebra}
\label {Appendix}

The proof of Lemma~\ref{matching}  required
{\em   finding $\hat g\in \hat L = \hat E_7(q)$ such that 
${\big (( S\cap  L)\Psi_L^{-1}\big )^{\hat g}}
\break
 = \hat S\cap \hat  L$.}
Since  $\hat A:= ( S\cap  L)\Psi_L^{-1} $ and $ \hat S\cap \hat  L$ are conjugate in
$\hat L$, we can use the behavior of the latter group in order to deduce properties of the
former one.
 
The group  $\cent_{\hat G}(\hat S)=\hat E_6(q)$ acts on the Lie algebra ${\mathcal L} \big(\hat E_7(q) \big)$
of
$\hat L$, decomposing it  as
$133=78\oplus 27
\oplus 27^*\oplus 1$, where
78 is   the Lie algebra ${\mathcal L} \big(\hat E_6(q) \big)$ of $\cent_{\hat G}(\hat S)$, the  27s are the usual dual pair of irreducible  
$\cent_{\hat G}(\hat S)$-modules of that dimension, 
and the 1-space is centralized.
The torus $\hat S\cap \hat  L$  centralizes $\cent_{\hat G}(\hat S)=\hat  E_6(q)$; each of its elements  acts
as a scalar $\rho$  on $27$  and  $\rho^{-1}$ on its dual  $ 27^*$,
so that $\hat S\cap \hat  L$ is nontrivial on both of those subspaces
(since $q>2$); and each of its elements   is 
 1 on  ${\mathcal L}\big(\hat E_6(q)\big)$
since each   is  both  an automorphism   of 
that algebra  and a scalar by Schur's
Lemma.  
Then $\hat S\cap \hat  L$ centralizes 
  $78\oplus 1$, so that
   $78$ is  the derived
Lie algebra  $
\cent_{{\mathcal L}(\hat E_7(q))} (\hat   S\cap \hat  L)'\cong {\mathcal L}\big(\hat E_6(q)\big)$.

With this background we proceed as follows. \!\! Find 
$ \cent_{{\mathcal L}(\hat E_7(q))} (\hat A) $
and then   $\cent_{{\mathcal L}(\hat E_7(q))} (\hat A) '\cong {\mathcal
L}\big(\hat E_6(q)\big)$, using elementary linear algebra.

Find a Chevalley basis $\{e_\alpha , e_{-\alpha}, h_\alpha \mid \alpha\in \Phi_6 \}$ of
$\cent_{{\mathcal L}(\hat E_7(q))} (\hat A) ' $ using
\cite{CM,CR}.  

Let $\Delta_6$  be a base for $\Phi_6  $.

Find the linear transformations $E_{ \alpha}(t)=\ad te_{-\alpha}$ and $
E_{-\alpha}(t)=\ad te_{-\alpha} $  for $\alpha\in \Delta_6$ and $t=f_k$ or $-f_k^{-1}$  in
$\F$; and then also $h_\alpha (f_k)$   as in (\ref{Definition of $h alpha$}).
Then  $\langle h_\alpha (\F^*) \mid \alpha \in \Delta_6 \rangle$ is a maximal split torus
of a group (isomorphic to $ \hat  E_6(q)$) of automorphisms of $\cent_{{\mathcal L}(\hat E_7(q))}
(\hat A) '$.

We saw above that  $\hat A$  is 1 on the $78$-space
$\cent_{{\mathcal L}(\hat E_7(q))} (\hat A) '$.   It follows that $T_7:=\langle h_\alpha (f_k) ,\hat A
\mid \alpha \in \Delta_6,  1\le k\le e \rangle$  is the direct product $
\langle h_\alpha (f_k)
\!\mid \!\alpha \in \Delta_6, {1\le k\le e}\rangle \times \langle \hat A\rangle $,
and hence  has the
correct order 
$(q-1)^6(q-1)$ to be a maximal torus of $\hat L$.

We can now  obtain a Chevalley basis of ${\mathcal L}\big(\hat E_7(q)\big)$: 
diagonalize the action of $T_7$ on ${\mathcal L}(\hat E_7(q))$ and normalize the basis as in  \cite[Sec.~4.2]{Ca1}.

We now have two Chevalley bases of ${\mathcal L}(\hat E_7(q))$: the one we started with 
(which was  implicitly used to write the generators of $\hat  E_7(q)$ in the
presentation  (\ref{relation1})--(\ref{relation4})), and the one just constructed.  Let
$\hat g$ be the linear transformation effecting the corresponding base change.  It is in
$\hat E_7(q)$, so we can use our
$E_7(q)$ algorithm for Theorem~\ref{Main Theorem}(iv)~  (a recursive call) to write it
using a straight-line program in the generators of
$\hat L$, as required.

\end{document}